\documentclass[twoside, 11pt]{article}

\usepackage{amssymb, amsmath, latexsym, mathrsfs, verbatim, calc, graphicx}

\numberwithin{equation}{section}
\newtheorem{thm}{Theorem}[section]
\newtheorem{lem}[thm]{Lemma}
\newtheorem{cor}[thm]{Corollary}
\newtheorem{prop}[thm]{Proposition}
\newtheorem{rem}[thm]{Remark}

\newtheorem{eg}[thm]{Example}
\newtheorem{defn}[thm]{Definition}

\def\ba{\begin{array}}
\def\ea{\end{array}}
\def\beq{\begin{equation}}
\def\eeq{\end{equation}}
\def\bes{\begin{equation*}}
\def\ees{\end{equation*}}
\def\bea{\begin{eqnarray}}
\def\eea{\end{eqnarray}}
\def\beaa{\begin{eqnarray*}}
\def\eeaa{\end{eqnarray*}}
\def\lan{\mathop{\langle}}
\def\ran{\mathop{\rangle}}

\def\dis{\displaystyle}

\def\no{\noindent}

\def\lastline{\par \vspace{-7.3ex} \no}

\def\nts{\negthinspace}

\def\ms{\medskip}
\def\bs{\bigskip}
\def\q{\quad}
\def\qq{\qquad}

\def\ua{\mathop{\uparrow}}
\def\da{\mathop{\downarrow}}

\def\={=\nts \nts=\nts \nts=\nts \nts=}
\def\cds{\cdots}

\def\wt{\widetilde}

\def\({\textnormal{(}}
\def\){\textnormal{)}}

\def\cd{\cdot}
\def\cds{\cdots}

\def\pa{\partial}

\def\a{\alpha}
\def\g{\gamma}

\def\z{\zeta}

\def\si{\sigma}
\def\vsi{\varsigma}
\def\t{\tau}
\def\f{\varphi}
\def\th{\theta}
\def\o{\omega}
\def\h{\widehat}
\def\f{\phi}

\def\G{\Gamma}

\def\O{\Omega}

\def\Th{\Theta}

\def\bF{{\bf F}}

\def\ba{{\bf a}}

\def\bx{{\bf x}}

\def\cA{{\cal A}}
\def\cB{{\cal B}}

\def\cF{{\cal F}}
\def\cG{{\cal G}}

\def\cM{{\cal M}}

\def\dbE{\rm l\nts E}

\def\dbN{\rm l\nts N}

\def\dbR{\rm I\nts R}

\def\hB{\mathbb{B}}

\def\hE{\mathbb{E}}

\def\hN{\mathbb{N}}

\def\hR{\mathbb{R}}

\def\neg{\negthinspace}

\def\tr{{\hbox{tr}}}

\def\qed{\hfill \rule[0cm]{.25cm}{.25cm}\medskip}   
\def\dfnn{\stackrel{\triangle}{=}}

\def\b1{{\bf 1}}

\newenvironment{proof}{\no {\bf Proof.\;}}{$\qed$\vspace{.1in}}

\newenvironment{itm}{\vspace{-1ex}\begin{itemize}}{\end{itemize}}
\def\bi{\begin{itm}}
\def\ei{\end{itm}}

\def\equ_ind{\arabic{section}.\arabic{equation}}
\def\sec_ind{\arabic{section}}

\begin{document}

\title{\bf
Pathwise Taylor Expansions for It\^o Random Fields}

\author{
Rainer Buckdahn\thanks{ \no D\'{e}partement de Math\'{e}matiques,
Universit\'{e} de Bretagne-Occidentale, F-29285 Brest Cedex, France;
email: Rainer.Buckdahn@univ-brest.fr. }, ~ Ingo Bulla\thanks{\no
Department of Bioinformatics, University of G\"{o}ttingen, 37077 G\"{o}ttingen; email:
ibulla@uni-goettingen.de. This author was supported by a fellowship of the DAAD
(German Academic Exchange Service).}, ~ and ~ Jin Ma\thanks{ \noindent Department of
Mathematics, University of Southern California, Los Angeles, 90089;
Department of Mathematics, Purdue University, West Lafayette, IN
47907-1395; email: jinma@usc.edu. This author is supported in part
by US NSF grant \#0505427 and \#0806017. Part of this work was
completed while this author was visiting the Department of
Mathematics, University of Bretagne-Occidentale, France, whose
hospitality was greatly appreciated.} }




\date{}
\maketitle

\begin{abstract}
In this paper we study the  {\it pathwise stochastic Taylor expansion}, in
the sense of our previous work \cite{Buckdahn_Ma_02}, for a class of
It\^{o}-type random fields in which the diffusion part is allowed to contain
both the random field itself and its spatial derivatives.
Random fields of such an ``self-exciting" type particularly contains the
fully nonlinear stochastic PDEs of curvature driven diffusion, as well as
certain stochastic Hamilton-Jacobi-Bellman equations. We introduce the new notion
of ``$n$-fold" derivatives of a random field, as a fundamental device to cope
with the special self-exciting nature. Unlike our previous work \cite{Buckdahn_Ma_02},
our new expansion can be defined around any random time-space point $(\t,\xi)$, where
the temporal component $\t$ does not even have to be a stopping time. Moreover, the exceptional
null set is independent of the choice of the random point
$(\t,\xi)$. As an application, we show how this new form of pathwise Taylor expansion
could lead to a different treatment of the stochastic characteristics for a class of fully
nonlinear SPDEs whose diffusion term involves both the solution and its gradient, and hence
lead to a definition of the  {\it stochastic viscosity solution} for such SPDEs, which is new
in the literature.

\end{abstract}

\vfill \bs

\no


\bs

\no{\it 2000 AMS Mathematics subject classification:} 60H07,15,30;
35R60, 34F05.



\section{Introduction}
\setcounter{equation}{0}\label{Introduction}

In our previous work \cite{Buckdahn_Ma_02} we studied the so-called
{\it pathwise stochastic Taylor expansion} for a class of
It\^{o}-type random fields. The main result can be briefly described
as follows. Suppose that $\{u(t,x), (t,x)\in[0,T] \times\hR^n\}$ is
an It\^{o}-type random field of the form
 \beq
 \label{itofield}
 u(t,x)=u(0,x)+\int_0^t
 u_1(s,x)ds+\int_0^t u_2(s,x)dB_s, \q (t,x)\in[0,T] \times\hR^n,
 \end{equation}
where $B$ is a 1-dimensional standard Brownian motion, defined on a
complete probability space $(\O,\cF,P)$. If we denote
$\bF=\{\cF_t\}_{t\ge0}$ to be the natural filtration generated by
$B$ and augmented by all $P$-null sets in $\cF$,  then under
reasonable regularity assumptions on the integrands $u_1$ and $u_2$,
 the following stochastic ``Taylor expansion" holds:
{\it For any stopping time $\t$ and any $\cF_\t$-measurable,
square-integrable random variable $\xi$, and for any sequence of
random variables $\{(\t_k,\xi_k)\}$ where $\t_k$'s are stopping
times such that either $\t_k>\t$, $\t_k\da \t$; or $\t_k<\t$,
$\t_k\ua \t$,
and $\xi_k$'s are all $\cF_{\t_k\wedge \t}$-measurable, square
integrable random variables, converging to $\xi$ in $L^2$, it holds
almost surely that
 \bea
 \label{itotaylor}
 u(\t_k,\xi_k)& =&u(\t,\xi)+a(\t_k-\t)+b(\xi_k-\xi)+{c\over
 2}(B_{\t_k} -B_\t)^2 \\
 &&\q+\lan p,\xi_k-\xi\ran+\lan q, \xi_k-\xi\ran(B_{\t_k}-B_\t)+\frac12 \lan
 X(\xi_k-\xi),\xi_k-\xi\ran \nonumber\\
 & &\q+o(|\t_k-\t|)+o(|\xi_k-\xi|^2), \nonumber
 \eea
where $(a,b,c,p,q,X)$ are all $\cF_\t$-measurable random variables,
and the remainder $o(\z_k)$ are such that $o(\z_k)/{\z_k} \to 0$ as
$k\to\infty$, in probability. Furthermore, the six-tuple
$(a,b,c,p,q,X)$ can be determined explicitly in terms of $u_1$,
$u_2$ and their derivatives.
}

By choosing $u_1$ and $u_2$ in different forms, we then extended the
Taylor expansion to solutions of stochastic differential equations
with initial state as parameters, and to solutions of nonlinear
stochastic PDEs. In the latter case we further introduced the notion
of the {\it stochastic super-(sub-)jets} using the Taylor expansion,
from which the definition of stochastic viscosity solution was
produced. We should note that in \cite{Buckdahn_Ma_02} all the SDEs
and SPDEs have the diffusion coefficient in the form
$g(t,x,u(t,x))$, that is, they only involve the solutions
themselves. Such a structure turns out to be essential for the
so-called Doss-Sussmann transformation, and in that case the
stochastic viscosity solution became natural.

In this paper we are interested in the stochastic Taylor expansion
for random fields of the following form:
 \bea
 \label{itofield2}
 u(t,x)=u(0,x)+\int_0^tu_1(s,x)ds+\int_0^t\lan g(t,x,Du(s,x)),dB_s\ran,
 \q (t,x)\in[0,T]\times\hR^n,
 \eea
where $u_1$ is a random field, and $g$ is a deterministic function. We shall
assume that they are ``smooth" in the sense that all the desired derivatives
exist, and the degree of smoothness will be specified later.

Given such a random field, we again consider
the possibility to expand $u$ in the sense of (\ref{itotaylor}), but this
time in a more natural way: We shall allow the pair
$(\tau,\xi):\Omega\rightarrow [0,T]\times \hR^d$ to  be arbitrary
random points.
Furthermore, we should note that although the Taylor expansion
(\ref{itotaylor}) holds almost surely, in general the null set may depend on
the choice of random point $(\tau,\xi)$. In this paper we shall look for
a universal expansion, in the sense that there is a subset $\widetilde{\Omega}
\subset \Omega$ with full probability measure, on which the stochastic
Taylor expansion holds for all choices of
random points $(\tau,\xi)$. These technical improvements make the stochastic
Taylor expansions much more ``user-friendly", and more importantly, it
will be more effectively used in our study of stochastic viscosity solution,
especially in the proof of uniqueness, as we shall see in our forthcoming
publications on that subject.

As one can easily observe from (\ref{itofield2}) that the random
field $u$ is actually already a solution to a first order stochastic
PDE. An immediate difficulty in deriving the Taylor expansion is the
characterization of the derivatives of the random field, as it will
increase in a ``bootstrap" way, very similar to those that one has
often seen in the anticipating calculus. As a consequence, the
original Doss-Sussmann type transformation used in our previous
works \cite{Buckdahn_Ma_02,Buckdahn_Ma_01_1,Buckdahn_Ma_01_2} no
longer works in this case. To overcome this difficulty we introduce
the notion of ''{\it $n$-fold derivative}" of a random field, which
essentially takes $(u,Du,\cds,D^{n-1}u)$ as a vector-valued random
field, and define its derivative in a recursive way. Such a
definition turns out to be very close to the idea of converting a
higher order ordinary differential equation to a first order system,
and is mainly motivated by the ``{\it stochastic characteristics}"
of a stochastic PDE (cf. e.g., Kunita \cite{Kunita_90}). In fact, by
combining the definition of stochastic characteristics in
\cite{Kunita_90} and the stochastic Taylor expansion developed in
this paper, we are able to rigorously define a stochastic
diffeomorphism that relates the Stochastic PDE of the form
\begin{equation}
    \label{Def. u}
\left\{\begin{array}{lll}
\dis du(t,x)   =
f(x,u,D_xu,D^2_{x}u)dt +  g(x, D_xu))\circ dB^i_t,
\quad (t,x)\in(0,T)\times\hR^d; \\
  u(0,x)     =  u_0(x), \quad x\in\hR^d,
\end{array}\right.
\end{equation}
to a PDE without Brownian components in their Taylor expansions
around any temporal-spatial point, generalizing the Doss-Sussmann
transformation in our previous works to the present case.



We should note that one of the main difficulties in the study of the
stochastic viscosity solution can be described as ``local" vs.
``global". That is, the local nature of the viscosity solution vs.
the global nature of the stochastic analysis (e.g., stochastic
integrals). Our idea is to ``localize" the stochastic integral, or
solution to the SPDE via the stochastic Taylor expansion that we
established in the previous section. It should be noted that the
universal set ``$\O'$" that we found in these expansions is the
essential point here.

Finally, we would like to point out that the topic of stochastic Taylor
expansion has been explored by many authors in various forms, and used to
provide numerical and other approximation schemes for SDEs and SPDEs or
randomized ODEs and PDEs (see, for example, \cite{azencott, arous,jentzen,platen},
to mention a few). These expansions often use either the Lie-algebraic structure
of  the path space or the chaotic type expansion of multiple stochastic integrals.
As a consequence it is hard to deduce the simultaneous spatial-temporal expansions
that we are pursuing in this work, especially when the remainders are
estimated in a pathwise manner. We should also note that the pathwise
version of stochastic viscosity solutions, suggested by Lions-Souganidis
\cite{Lions_Souganidis_98_1, Lions_Souganidis_98_2}, has found an effective
framework recently, using the theory of {\it rough path} (cf. Caruana-Friz-Oberhauser \cite{CFO}).
However, due to the special nature of the rough path integrals, the arguments seem to depend 
heavily on the fact that there exist stochastic characteristics in the form of 
$C^3$-diffeomorphisms which transform the SPDE to a pathwise PDE. Consequently, the SPDE 
studied in \cite{CFO}, while fully nonlinear in the drift, seems to be restricted to the cases
when the diffusion coefficient depends only (linearly) on $Du$, the gradient of the solution, so that a chain-rule type of argument could be applied.
The generality of the diffusion part in the fully
non-linear PDE suggested in this paper, and the stochastic characteristics related to it,
does not seem to be an easy consequence of such a method. Furthermore, our Taylor
expansions are constructed within a more ``elementary" stochastic analysis framework, by exploiting the
properties of Brownian motion without using the advanced algebraic geometric structure of the path spaces, therefore we believe that it provides a more accessible alternative.

This paper is organized as follows. In Section \ref{Preliminaries
and notations} we clarify all the necessary notations, and
state the main theorem. In section 3 we give a fundamental estimate
of this paper, regarding multiple stochastic integrals. In sections
4 and 5 we study the forward and backward Taylor expansion,
respectively. We note that in the stochastic case, the temporal
direction of the expansion does affect the outcome. Finally in
section 6 we try to apply the Taylor expansion to the solution of a
class of nonlinear SPDEs, and propose a possible definition of the
stochastic viscosity solution in this case.


\section{Preliminaries and Statement of the Main Theorem}
\label{Preliminaries and notations}
\setcounter{equation}{0}

Throughout this paper we denote $(\Omega, \cF,P)$ to be a complete
probability space on which is defined an $\ell$-dimensional Brownian
motion $B=(B_t)_{t\ge 0}$. Let $\bF^B\dfnn\{\cF^B_t\}_{t \ge 0}$ be
the natural filtration generated by $B$, augmented by the $P$-null
sets of $\cF$; and let $\cF^B =\cF^B_\infty$. We denote
$\cM_{0,T}^B$ to be the set of all  $\bF^B$-stopping times $\tau$
such that $0\leq\tau\leq T$, $P$-a.s., where $T>0$ is some fixed
time horizon; and denote $\cM_{0,\infty}^B$ to be all
$\bF^B$-stopping times that are almost surely finite.

In what follows we write $\hE$ (also $\hE_1$, $\cds$) for a generic
Euclidean space, whose inner products and norms will be denoted as
the same ones $\langle\cdot,\cdot\rangle$ and $|\cdot|$,
respectively; and write $\hB$ for a generic Banach space with norm
$\|\cd\|$. Moreover, we shall denote $\cG\subseteq\cF^B$ to be a
sub-$\si$-field of $\cF^B$, and for any $x\in\hR^d$ and constant
$r>0$ we denote $\overline{B}_r(x)$ to be the
closed ball with center $x$ and radius $r$.
 Furthermore, the following spaces of functions will be frequently
used in the sequel. We denote

\begin{itemize}

\item $L^p(\cG;\hE)$ to be all $\hE$-valued,
$\cG$-measurable random variables $\xi$, with $E(|\xi|^p\}<\infty$.
Further, we shall denote
$L^{\infty}_-({\cal G};\hE)\dfnn\cap_{p>1}
L^p(\cG;\hE)$.


\item $L^q(\bF^B,[0,T];\hB)$ to be all
$\hB$-valued, $\bF^B$-progressively measurable processes $\psi$,
such that $E\int_0^T\|\psi_t\|^qdt<\infty$. In particular, $q=0$
stands for all $\hB$-valued, $\bF^B$-progressively measurable
processes; and $q=\infty$ denotes all processes in
$L^0(\bF^B,[0,T];\hB)$ that are uniformly bounded.

\item $C^{k,\ell}([0,T]\times\hE;\hE_1)$ to be the space
of all $\hE_1$-valued functions defined on $[0,T]\times\hE$ which
are $k$-times continuously differentiable in $t\in[0,T]$ and
$\ell$-times continuously differentiable in $x\in\hE$.

\item $C^{k,\ell}_b([0,T]\times\hE;\hE_1)$,
$C^{k,\ell}_l([0,T]\times\hE;\hE_1)$,
and $C^{k,\ell}_{p}([0,T]\times\hE;\hE_1)$, etc.
to be the subspace of $C^{k,\ell}([0,T]\times\hE;\hE_1)$, where the
subscript ``$b$" means all functions and their partial derivatives
are uniformly bounded; ``$l$" means
all functions are of at most linear growth; and
``$p$" means all functions and their partial derivatives are of at
most polynomial growth. The subspaces with the combined subscripts
of $b$, $l$, and $p$ are defined in an obvious way.

\item for any sub-$\sigma$-field $\cG\subseteq \cF^B_T$,
$C^{k,\ell}(\cG,[0,T]\times\dbE;\dbE_1)$ (resp.
$C^{k,\ell}_b(\cG,[0,T]\times\dbE;\dbE_1)$,
$C^{k,\ell}_p(\cG,[0,T]\times\dbE;\dbE_1)$) to be the space of all
$C^{k,\ell}([0,T]\times\dbE;\dbE_1)$ (resp.
$C^{k,\ell}_b([0,T]\times\dbE;\dbE_1)$, $C^{k,\ell}_p(
[0,T]\times\dbE;\dbE_1)$)-valued random variables that are
${\cG}\otimes\cB([0,T]\times \dbE)$-measurable;

\item $C^{k,\ell}(\bF^B,[0,T]\times\hE;\hE_1)$ (resp.
$C^{k,\ell}_b(\bF^B,[0,T]\times\hE;\hE_1)$, $C^{k,\ell}_p(\bF^B,
[0,T]\times\hE;\hE_1)$) to be the space of all random fields
$\varphi\in C^{k,\ell}(\cF^B_T,[0,T]\times\hE;\hE_1)$ (resp.
$C^{k,\ell}_b(\cF^B_T, [0,T]\times\hE;\hE_1)$,
$C^{k,\ell}_p(\cF^B_T,[0,T]\times\hE;\hE_1)$), such that for fixed
$x\in\hE$, the mapping $(t,\omega)\mapsto \varphi(t,x,\omega)$ is
$\bF^B$-progressively measurable, and for $P$-a.e. $\o$, the mapping
$\f(\cd,\cd,\o)\in C^{k,\ell}([0,T]\times\hE;\hE_1)$ (resp.
$C^{k,\ell}_b([0,T]\times\hE;\hE_1)$, $C^{k,\ell}_p([0,T]\times\hE;\hE_1)$).


\end{itemize}

If $\hE_1=\hR$, we shall drop $\hE_1$ from the notation (e.g.,
$C^{k,\ell}([0,T]\times\hE)$, and so on); and we write
$C^{0,0}([0,T]\times\hE;\hE_1)= C([0,T]\times\hE;\hE_1)$, and
$C^{0,0}(\bF^B,[0,T]\times\hE)=C(\bF^B,[0,T] \times\hE)$,..., etc.,
to simplify notation. Finally, for $(t,x,y)\in [0,T]
\times\hR^d\times\hR$, we denote $D = D_x =\nabla_x= \bigl(
\frac{\partial}{\partial x_1},\cdots, \frac{\partial}{\partial x_d}
\bigr)$, $D^2=D_{xx}=(\partial^2_{x_i x_j})_{i,j=1}^d$,
$D_y=\frac{\partial}{\partial y}$, and $D_t=\frac{\partial}{\partial
t}$. The meaning of $D_{xy}$, $D_{yy}$, etc., should be clear.

%

Finally, since the random fields that we are interested in are
always of the form of (\ref{itofield2}), which is an SPDE already,
the following definition in \cite{Buckdahn_Ma_02} is useful.
Consider the fully nonlinear second-order SPDE:
\begin{equation}
\label{exspde} u(t,x)=u_0(x)+\int_0^tf(s,
x,(u,Du,D^2u)(s,x))ds+\int_0^tg(t,x,(u,Du)(s,x))dB_s,
\end{equation}
where $(t,x)\in[0,T]\times \hR^d$, and $f$ and $g$ are functions
with appropriate dimensions.

\begin{defn} A random field $u=\{u(t,x,\o): (t,x,\o)\in[0,T]
\times\hR^d\times\O\}$ is called a ``regular'' solution to SPDE
(\ref{exspde}) if

(i) $u\in C^{0,2}(\bF^B;[0,T]\times \hR^d)$;

(ii) $u$ is an It\^{o}-type random field with the form
$$u(t,x)=u_0(x)+\int_0^tu_1(s,x)ds+\int_0^t u_2(s,x)dB_s,\q
(t,x)\in[0,T]\times \hR^d,
$$
where
$$ u_1(t,x)=f(t,x,(u,Du,D^2u)(t,x)),\q u_2(t,x)=g(t,x,(u,Du)(t,x)),
$$
for all $(t,x)\in[0,T]\times\hR^d$, $P$-a.s. \qed
\end{defn}


In this paper we will consider a special type of ``smoothness" for a
random field, defined through what we shall call the ``$n$-fold
differentiability" below. Such a characterization of
differentiability is mainly motivated by the ``{\it stochastic
characteristics}" for stochastic PDEs (cf. e.g., \cite{Kunita_90}),
which often take the form of a system of first order Stochastic
PDEs. The idea is actually quite similar to the well-known
transformation from a higher order ordinary differential equation
(ODE) to a first order system of ODEs.

To begin with, we recall that for any multi-index $j=(j_1,\cds,j_d)$, its ``length"
is defined by $|j|\dfnn \sum_{k=1}^d j_k$. We have the following definition.
\begin{defn}
\label{nfold} A random field $\zeta\in C^{0,n}(\bF^B,[0,T]\times
\hR^d)$ is called ``$n$-fold" differentiable in the spatial variable
$x$ if there exist
$n$ smooth random fields $\zeta_i\in C^{0,n}(\bF^B,[0,T]\times
\hR^d;\hR^{d_i})$, $2\leq i\leq n+1$, with $d_1=d$ and $d_i\in\hN$,
$2\leq i\leq n+1$, and the functions $F_i,G_i:\hR^d\times
\hR^{d_{i+1}}\rightarrow \hR^d$, $i=1,\cds, n$,
such that, denoting
$\zeta_1\dfnn\zeta$, the following properties are satisfied:
\begin{itemize}
\item[(T1)] $F_i,G_i\in C^\infty_{\ell,
p}$, $i=1,\cds,n$;

\item[(T2)]  For $1\leq i\leq n$, it holds that
 \bea
 \label{zetai}
  \zeta_i(t,y)=\zeta_{i,0}(y)+\int_0^{t}F_{i}(y,\zeta_{i+1}(s,y))ds
 +\int_0^{t}G_{i}(y,\zeta_{i+1}(s,y))dB_s,
 \eea
for all $(t,y)\in [0,T]\times \hR^d$, with $\zeta_{i,0}\in
C^2(\hR^d;\hR^{d_i})$, $1\leq i\leq n$.

\item[(T3)]
For any $1\leq i\leq n+1$, $m\geq 1$, and multi-index
$j=(j_1,...,j_d)$ with $\vert j\vert=j_1+...+j_d\leq n$, it holds
that
$$\sup\{\vert D_y^j\zeta_i(t,y)\vert, t\in
[0,T], \vert y\vert\leq m\}\in L^{\infty}_-(\Omega, {\cal F},P).
$$
\end{itemize}
We shall call $\z_i$, $i=2,\cds,n+1$ the ``generalized derivatives"
of $\z=\z_1$, with ``coefficients" $(F_i,G_i)$, $i=1,\cds,n$.
\end{defn}

For notational convenience
we will often write $F=F_1$ and $G=G_1$ when there is no danger of
confusion. We denote the set of all $n$-fold differentiable random
fields by $C^{0,(n)}_\bF([0,T]\times \hR^d)$.


Before we state the main result of this paper, let us note again
that the main feature of our stochastic Taylor expansion is that it
can be expanded around any random point $(\t,\xi)$, and that the
expansion holds almost surely with an exceptional set that is
independent of the choice of $(\t,\xi)$. But this amounts to
saying that the expansion can be performed around any
deterministic point $(\t,\xi)$ with any (deterministic)
increments, outside a uniform exceptional set.
In other words, the complicated approximating sequences
$(\t_k,\xi_k)$ in (\ref{itotaylor}) can be replaced by simple
deterministic increments $(t+h, x+k)$, for all $(h,k)$ near zero. We
should also note that for the purpose of our application, in this
paper we consider only the Taylor expansion up to the second order,
and for that purpose the 3-fold differentiability of the random
field would suffice. The precise statement of our main result is the
following theorem.



\begin{thm}
\label{StoTaylor} Let $\zeta\in C^{0,(3)}(\bF^B, [0,T]\times \hR^d)$
be a random field satisfying the standard assumptions (T1)-(T3) with
generalized derivatives $\z_i$, $i=2,3,4$ and coefficients
$(F_i,G_i)$, $i=1,2,3$. Then, for every $\alpha \in
\big(\frac{1}{3},\frac{1}{2}\big)$ and every $m\in\hN$ there exist a
subset $\tilde{\O}_{\a,m}\subset\Omega$ with
$P\{\tilde\O_{\a,m}\}=1$, such that, for all $(t,y,\o)\in
[0,T]\times \overline{B}_m(0)\times\tilde\O_{\a,m}$,
we have the following Taylor expansion
 \bea
 \label{taylorexp}
 \zeta(t+h,y+k) -
 \zeta(t,y)  \neg\neg&\neg =\neg &\neg\neg ah + b(B_{t+h} - B_t) + \frac{c}{2} (B_{t+h} -
 B_t)^2 + \langle p,k \rangle
      + \frac{1}{2}\langle Xk,k \rangle \\
 &   & + \langle q,k \rangle (B_{t+h} - B_t) + (|h| +
 |k|^2)^{3\alpha}R_{\alpha,m}(t,t+h,y,y+k), \nonumber
\eea
for all $(t+h, y+k)\in [0,T]\times \overline{B}_m(0)$. Here,
with $(F,G)=(F_1,G_1)$, one has
 \bea
 \label{taylors}
  a & = & F(y,\zeta_2(t,y)) - \frac{1}{2} (D_z G)(y,\zeta_2(t,y))G_2(y,\zeta_3(t,y)),
  \nonumber\\
  b & = & G(y,\zeta_2(t,y)), \q
  c  =  (D_z G)(y,\zeta_2(t,y))G_2(y,\zeta_3(t,y)), \nonumber\\
  p & = & D_y\zeta(t,y), \q
  X  =  D_y^2\zeta(t,y), \\
  q & = & (D_y G)(y,\zeta_2(t,y)) + (D_z
  G)(y,\zeta_2(t,y))D_y\zeta_2(t,y). \nonumber
 \eea
Furthermore, the remainder of Taylor expansion $R_{\alpha,m}:
[0,T]^2\times (\hR^d)^2\times\Omega\mapsto \hR$ is a measurable
random field such that
 \begin{eqnarray}
 \label{Ram}
 \overline{R}_{\a,m}\dfnn\sup_{t,s\in[0,T];~
 y,z\in\overline{B}_m(0)}|R_{\alpha,m}(t,s,y,z)|
 \in L^{\infty}_{-}(\O, \cF, P).
 \end{eqnarray}
\end{thm}

{\it Proof.} Since the proof of this theorem is quite lengthy and
technical, we shall split it into several cases and carry it out in
the following sections. \qed

\begin{rem}
\label{randompt} {\rm (i) It is worth noting that the ``universal"
estimate for the remainder $R_{\a,m}$ is the main reason why the
Taylor expansion can now hold around any random space-time point. In
other words, the points $(t,y)$ and $(h, k)$ in Theorem
\ref{StoTaylor} can be replaced by any $(\t,\xi)$, $(\si,\eta)\in
L^0(\cF^B;[0,T])\times L^0(\cF^B;\hR^d)$, such that
$(\t+\si,\xi+\eta)\in[0,T]\times \overline{B}_m(0)$, $P$-a.s. on
$\tilde\O$, except in that case the reminder should read
$$ \hat
R_{\a,m}\dfnn R_{\a,m}(\t,\t+\si,\xi,\xi+\eta), \q \mbox{where} \q
|\hat R_{\a,m}|\le \overline{R}_{\a,m}\in L^{\infty}_-(\bF^B;P).
$$
 In what follows
 we denote $R_{\alpha,m}$ to be a generic term satisfying
(\ref{Ram}), which is allowed to
vary from line to line.

(ii) From the expressions (\ref{taylors}) it is clear that the drift
term $F (=F_1)$ appears only in the coefficient $a$, and it is in
its original form. In fact, as we shall see in the proof, the exact
form of $F(t,y)(=F(y,\zeta_2(t,y)))$ is not important at all. We
could simply change it to $F(t,y)$ and the results remains the same.

(iii) Theorem \ref{StoTaylor} remains true if the assumption (T1) is
replaced by the weaker assumption:
\begin{itemize}
\item[(T1')]  $F_j\in C^{5-j}_{\ell, p}$;  $G_j\in
C^{7-j}_{\ell, p}$;  $1\le j\le 3$.
\end{itemize}
However, we shall not pursue this generality in this paper due to
the length of the paper.
\qed}
\end{rem}

To end this section we give an example, which more or less
motivated our study.
\begin{eg}
\label{egspde} {\rm Let $u=\{u(t,x)\}$ be a regular solution of the
SPDE (\ref{exspde}), and for simplicity we assume that both $f$ and
$g$ are ``time-homogeneous" (i.e., they are independent of the
variable $t$), and $g$ is independent of $x$ as well. We define
$\zeta=\zeta_1=u$, and $\zeta_{i+1}=(\zeta_i,D\zeta_i,D^2\zeta_i)$,
$1\le i\le 3$. Then, assuming that the coefficients $f$ and $g$ are
sufficiently smooth and their derivatives of all order are bounded,
one can show that (T1)--(T3) are satisfied.
%
Furthermore, applying Theorem \ref{StoTaylor} we see that on some
subset $\widetilde\Omega$ of  $\Omega$ with $P(\widetilde\Omega)=1$,
independent of the expansion point $(t,x)\in[0,T]\times \hR^d$, the
stochastic Taylor expansion (\ref{taylorexp}) holds for $u$, with
\begin{eqnarray}
\label{SPDEexp}
 \nonumber a& = &
 f(x,(u,Du,Du^2)(t,x)) - \frac{1}{2}\Big\{ (gD_ug)((u,Du)(t,x)) \\
\nonumber &   &  + D_ug((u,Du)(t,x))
      \langle D_pg((u,Du)(t,x)),Du(t,x) \rangle \\
\nonumber &   &  + \langle D^2u(t,x)D_pg((u,Du)(t,x)),
D_pg((u,Du)(t,x))\rangle
      \Big\} \nonumber \\
b& = & g((u,Du)(t,x))
\nonumber\\
c&  = &  (gD_ug)((u,Du)(t,x))  + D_ug((u,Du)(t,x))
      \langle D_pg((u,Du)(t,x)),Du(t,x)\rangle \\
\nonumber &   &+ \langle D^2u(t,x)D_pg((u,Du)(t,x),
        D_pg((u,Du)(t,x))\rangle
\nonumber\\
p & =  & Du(t,x); \q X= D^2u(t,x)
\nonumber\\
q & =  &  D_ug((u,Du)(t,x)Du(t,x)
      + D^2u(t,x) D_pg((u,Du)(t,x)), \nonumber
\end{eqnarray}
for all $h\in \hR$ and $k\in \hR^d$ such that $t+h\in [0,T]$ and
$|k|\le m$, respectively. Here, the coefficients $F_2, G_2, F_3$,
and $G_3$
 are obtained
by differentiating (\ref{exspde}) with respect to $x$.

The explicit expressions of the pathwise Taylor expansion in
(\ref{SPDEexp}) will be important for our study of the stochastic
characteristics for SPDE (\ref{exspde}), which will be discussed in
the last section of this paper.
}\qed
\end{eg}

\section{Some fundamental estimates}
\setcounter{equation}{0}

\label{Proof of the stochastic Taylor expansion}

Before we prove the theorem, let us first introduce the so-called
Wick-square of the Brownian motion, which is originated in the
Wiener homogeneous chaos expansion (cf., e.g., \cite{Hida},
\cite{Nua}).
For any $0\leq s\leq t$, we define the Wick-square of $B_t-B_s$ to
be
\[
  (B_t-B_s)^{\diamond 2} := (B_t-B_s)^2-|t-s|.
\]
Moreover, let $L_{\bF^B}^0(\Omega;
W_{2,loc}^{0,1}([0,T]\times\hR^d))$ denote the space of all random
fields $f\in L^0(\bF^B,[0,T]\times\hR^d)$ such that
$f(\omega,\cdot,\cdot)$ is $P$-a.s.~an element of
$$ W_{2,loc}^{0,1}([0,T]\times\hR^d)
   =  \{ u\in L^2_{loc}([0,T]\times\hR^d): \; D_x u\in
   L^2_{loc}([0,T]\times\hR^d)\},
   $$
where in this case $D_x u$ denotes the weak partial derivative of $u$ with respect to
 $x$.

We begin by proving an important estimate for multiple Stochastic
integrals. We consider the multiple stochastic integral defined
recursively as follows. Let $f_j\in L_{\bF^B}^0(\Omega;
W_{2,loc}^{0,1}([0,T]\times\hR^d))$, for $1\le j\le N$, $N\in \hN$.
Let $X^0_{t,s}(x)\equiv1$, and for $j=1,\cdots, N$,
 \bea
 \label{Xstj} X^j_{t,s}(x)=\int_t^s f_j(s_j,x)X^{j-1}_{t,s_j}(x)dB_{s_j}, \q s\in [t,T].
 \eea
It is also useful to define the similar multiple integrals: for
$k>0$, let $X^{k,k-1}_{t,s}(x)\equiv1$, and for $l\ge k$,
 \bea
 \label{Xstkl}
 X^{k,l}_{t,s}(x)=\int_t^s f_{l}(u_l,x)X^{k,l-1}_{t,u_l}(x)dB_{u_l}, \q s\in [t,T].
 \eea
Then it is clear that $X_{t,s}^{1,j}=X_{t,s}^j$, for $j\ge 1$. We
have the following {\it regularity} result.

\begin{prop}
\label{Estimation iterated Ito integral} Let $N\in\hN$ be given.
Assume that   $f_j\in L_{\bF^B}^0(\Omega;
W_{2,loc}^{0,1}([0,T]\times\hR^d))$,  $1 \leq j \leq N$, satisfy
that
 \bea
 \label{assmpf}
  C_{m,p}(f_j)\dfnn E\left[ \int_{|x|\leq m} \int_0^T |D_x^i f_j(t,x)|^p dt dx \right] <
  \infty,
 \eea
for all $m\in\hN$, $p\geq 1$, and $i=(i_1,\ldots,i_d)$ with
$|i|=i_1+\ldots+i_d \leq 1$. Then it holds that
\begin{equation}
\label{zeta_{beta,m} in L^infty-}
  \zeta_{\beta,m} \dfnn \sup\left\{ \frac{|X^N_{t,s}(x)|}{|s-t|^\beta}, \;
  0 \leq t < s \leq T, \; |x| \leq m \right\} \in L^{\infty}_-(\Omega, {\cal F}, P)
\end{equation}
for all $\beta\in\left(0,\frac{N}{2}\right)$, $m\in\dbN$.
\end{prop}


{\it Proof.}
The proof is based on the Kolmogorov continuity criterion for random
fields (cf. e.g., \cite[Theorem I.2.1]{Revuz_Yor_91}), combined with
an induction argument. We note that throughout the proof we shall
use the notations $C_m$, $C_{m,p}$, $C_{M,k}$, etc., to represent
the generic constants depending only on $f_i$'s, $T$,  and the
parameters in their subscripts, which are allowed to vary from line
to line.


We begin with the case $N=1$. In this case we write $f_1=f$ and
$X_1=X$ for simplicity. Let $m\in\hN$, $\gamma>2$, and $0\leq t < s
\leq T$. Applying Sobolev's embedding theorem, H\"older Inequality,
and noting that $X^0_{t,s}=1$, we have
\begin{eqnarray*}
&      & E\left[ \sup_{|x|\leq m}|X_{t,s}(x)|^\gamma \right] \\
& \leq & C_m E\left[ \left( \sum_{|j|\leq 1} \int_{|x|\leq m}
|D_x^jX_{t,s}(x)|^{d+1} dx
         \right)^{\frac{\gamma}{d+1}} \right] \\
& \leq & C_m \sum_{|j|\leq 1} \left( \int_{|x|\leq m} E\left[ \left|
\int_t^s   D_x^j f(r,x) dB_r
                 \right|^{\gamma(d+1)} \right] dx \right)^{\frac{1}{d+1}} \\
& \leq & C_m \sum_{|j|\leq 1} \left( \int_{|x|\leq m} E\left[ \left(
\int_t^s   |D_x^j f(r,x)|^2 dr
                 \right)^{\frac{\gamma(d+1)}{2}} \right] dx \right)^{\frac{1}{d+1}} \\
 & \leq & C_m \sum_{|j|\leq 1} \left( \int_{|x|\leq m}
E\left[ \left( (s-t)^\frac{\gamma(d+1) - 2}{\gamma(d+1)}
                 \left( \int_t^s|D_x^j f(r,x)|^{\gamma(d+1)} dr \right)^\frac{2}
                 {\gamma(d+1)} \right)^\frac{\gamma(d+1)}{2}
                 \right] dx \right)^{\frac{1}{d+1}}\\
 & \leq & C_m \sum_{|j|\leq 1} \left( \int_{|x|\leq m} E\left[
(s-t)^{\frac{\gamma(d+1)}{2} - 1}
                 \int_t^s|D_x^j f(r,x)|^{\gamma(d+1)} dr \right] dx \right)^{\frac{1}{d+1}} \\
& \leq & C_m (s-t)^{\frac{\gamma}{2}-\frac{1}{d+1}} \sum_{|j|\leq 1}
\left( E\left[ \int_{|x|\leq m}
                 \int_0^T|D_x^j f(r,x)|^{\gamma(d+1)} dr dx \right]
                 \right)^{\frac{1}{d+1}}.
\end{eqnarray*}
Consequently,
\[
  E\left[ \sup_{|x|\leq m}|X_{t,s}(x)|^\gamma \right] \leq C_{m,\gamma}(s-t)^{\frac{\gamma}{2}-\frac{1}{d+1}}.
\]
Endowing the space $C(\overline{B}_m)$ of continuous functions
defined on the closed $m$-ball $\overline{B}_m := \{ x\in\hR^d: \;
|x| \leq m \}$ with the sup-norm, and considering $\{x \mapsto
X_{0,s}(x); \, |x|\leq m\}_{s\in[0,T]}$ as a
$C(\overline{B}_m)$-valued process, we see that
\[
  \sup_{|x|\leq m}|X_{t,s}(x)| = |X_{0,s}(\cdot) -
  X_{0,t}(\cdot)|_{C(\overline{B}_m)}.
\]
Applying the Kolmogorov continuity criterion, we conclude that
$\zeta_{\beta,m}\in L^\gamma(\Omega, {\cal F}, P)$ for all $\beta
\in \left[ 0, \frac{1}{2}-\frac{d+2}{\gamma(d+1)} \right)$. Since we
can choose $\gamma>2$ arbitrarily large, (\ref{zeta_{beta,m} in
L^infty-}) follows.

We now prove the inductional step. That is, we assume that
(\ref{zeta_{beta,m} in L^infty-}) is true for $N-1$, and show that
it is also true for $N$. To do this, we shall adapt the proof of the
Kolmogorov continuity criterion given in \cite{Revuz_Yor_91} to our
framework. Let $0 \leq t < s \leq T$, $\gamma>2$, and $m\in\hN$.
First note that by a simple application of Burkholder-Davis-Gundy
and H\"older inequalities we have
 \bea
 \label{Xstest}
 E\Big\{\sup_{r\in[t,s]}|X_{t,r}^N(x)|^p\Big\}
 &=&E\Big\{\sup_{r\in[t,s]}
 \Big|\int_t^r f_N(s_N,x)X^{N-1}_{t,s_N}(x)dB_{s_N}\Big|^p\Big\}\nonumber\\
 &\le& C_{p}E\Big\{\Big[\int_t^s |f_N(s_N,x)|^2|X^{N-1}_{t,s_N}(x)|^2ds_N\Big]
 ^{p/2}\Big\}\\
  &\le& C_{p}E\Big\{\|f_N(\cd, x)\|^p_{L^2([t,s])}\sup_{s_N\in[t,s]}|
   X^{N-1}_{t,s_{N}}(x)|^{p}\Big\}
  \nonumber\\
  &\le& C_{p}\Big\{E\|f_N(\cd, x)\|^{pN}_{L^2([t,s])}\Big\}^{1/N}\Big\{
  E\Big\{\sup_{s_N\in[t,s]}|X^{N-1}_{t,s_{N}}(x)|^{\frac{pN}{N-1}}
  \Big\}^{\frac{N-1}{N}}.
  \nonumber
  \eea
By simply iterating the above argument we obtain that
 \bea
 \label{Xstest1}
 &&E\Big\{\sup_{r\in[t,s]}|X_{t,r}^N(x)|^p\Big\}^N \nonumber \\
 &\le & C_{p} E\|f_N(\cd, x)\|^{pN}_{L^2([t,s])}
 \cdot E\|f_{N-1}(\cd, x)\|^{pN}_{L^2([t,s])}\cdot
  E\Big\{\sup_{s_N\in[t,s]}\Big|X^{N-1}_{t,s_{N}}(x)|^{\frac{pN}{N-2}}
  \Big\}^{N-2}\\
 &\le& \cdots ~ \le C_{p}\displaystyle\, \prod_{k=1}^N\Big\{E\|f_k(\cd,
 x)\|^{pN}_{L^2([t,s])}\Big\}. \nonumber
 \eea
 Furthermore, with a similar argument as above we can also show that
for all $M,k\in\hN$, $p>1$, and multi-index $i$ satisfying $|i|=1$,
\begin{eqnarray*}
  && E\left\{ \sup_{t\leq r \leq s} |D^i X_{t,r}^{M}(x)|^k \right\} \\
&\le& E\left\{ \sum_{|i^1|,\ldots,|i^M|\leq 1} \sup_{t\leq r \leq s}
         \left| \int_t^r D^{i^M}f_M
         \Big(\int_t^{s_M}\cds
\Big(\int_t^{s_2}D^{i^1}f_1
dB_{s_1}\Big)\cdots dB_{s_{M-1}}\Big)dB_{s_M} \right|^k \right\} \\
& \leq & C_{M,k}\sum_{|i^1|,\ldots,|i^M|\leq 1} \prod_{j=1}^M
E\left[
          \|D^{i^j}f_j(s_j,x)\|_{L^2((t,s])}^{Mk}\right]^{\frac{1}{M}}
          \\
& \leq & C_{M,k}(s - t)^\frac{Mk(p-1)}{2p}\sum_{|i^1|,\ldots,|i^M|\leq 1} \prod_{j=1}^M
          E\left[\|D^{i^j}f_j(s_j,x))\|_{L^{2p}([t,s])}^{Mk}\right]^{\frac{1}{M}}.
\end{eqnarray*}
Here, we applied the H\"older inequality in the last step above.
Consequently, if $p\leq\frac{Mk}{2}$, then the assumption
(\ref{assmpf}) implies that, for all $m>0$,
\begin{equation}
\label{intEsupDX}
    \int_{|x|\leq m} E\left[ \sup_{t\leq r \leq s} |D^i X_{t,r}^{M}(x)|^k
    \right]dx
    \leq C_{m,M,k}(s - t)^\frac{Mk(p-1)}{2p}.
\end{equation}
We can then conclude that, for all $0 \leq t < s \leq T$, $\g>2$,
\begin{eqnarray*}
&      & E\left[ \sup_{|x|\leq m} |X_{t,s}^{N}(x)|^\gamma \right] \\
& \leq & C_m E\left[ \left( \sum_{|j|\leq 1} \int_{|x|\leq m} |D_x^j
X_{t,s}^{N}(x)|^{d+1} dx \right)^{\frac{\gamma}{d+1}}
                 \right] \\
& \leq & C_m E\left[ \left( \sum_{|i|,|j|\leq 1} \int_{|x|\leq m} \left| \int_t^s D_x^iX_{t,r}^{N-1}(x) D_x^j f_N(r,x)
         dB_r \right|^{\gamma(d+1)} dx \right)^{\frac{1}{d+1}} \right] \\
& \leq & C_m \left( \sum_{|i|,|j|\leq 1} \int_{|x|\leq m} E\left[ \left( \int_t^s |D_x^iX_{t,r}^{N-1}(x)|^2
         |D_x^j f_N(r,x)|^2 dr \right)^{\frac{\gamma(d+1)}{2}} \right] dx \right)^{\frac{1}{d+1}} \\
& \leq & C_m \left( \sum_{|i|,|j|\leq 1} \int_{|x|\leq m} E\left[
         \sup_{t\leq r\leq s} |D_x^iX_{t,r}^{N-1}(x)|^{\gamma(d+1)}
         \|D_x^j f_N(\cdot,x)\|_{L^2([t,s])}^{\g(d+1)} \right]
         dx \right)^{\frac{1}{d+1}}.
\end{eqnarray*}
Consequently, by H\"older inequality and (\ref{intEsupDX}) with
$p=\frac{Mk}{2}=(N-1)\gamma(d+1)$, we have
\begin{eqnarray*}
&      & E\left[ \sup_{|x|\leq m} |X_{t,s}^{N}(x)|^\gamma \right] \\
& \leq & C_m \left\{\sum_{|i|,|j|\leq 1} \Big[\int_{|x|\leq m}
E\Big[
         \sup_{t\leq r\leq s} |D_x^iX_{t,r}^{N-1}(x)|^{2\gamma(d+1)} \Big] dx \Big]^{\frac{1}{2}}
          \right. \times \\
&      & \phantom{C_m\left( \sum_{|i|,|j|\leq 1} \right.} \left.
\times \left( \int_{|x|\leq m}E\left[ \|D_x^j
         f_N(\cdot,x)\|_{L^2([t,s])}^{2\gamma(d+1)}
         \right]dx \right)^{\frac{1}{2}} \right\}^{\frac{1}{d+1}} \\
& \leq & C_{m,\gamma} \sum_{|i|\leq 1, \; j=1,2} \Bigg( (s-t)^{\frac12(N-1)\gamma(d+1) -
 \frac{1}{2}} \times \\
&      & \times\left( \int_{|x|\leq m} E\left[ \left( (s-t)^{\frac{\gamma(d+1)-1}{\gamma(d+1)}}
         \|D_x^i f_j(\cdot,x)\|_{L^{2\gamma(d+1)}([t,s])}^2\right)^{\gamma(d+1)} \right] \right)^\frac{1}{2} dx \Bigg)^{\frac{1}{d+1}} \\
& \leq & C_{m,\gamma} (s-t)^{\frac{N}{2}\gamma - \frac{1}{d+1}}
         \sum_{|i|\leq 1, \; j=1,2} \left( E\left[ \int_{|x|\leq m}
         \|D_x^i f_j(\cd,x)\|_{L^{2\gamma(d+1)}([t,s])}^{2\g(d+1)}dx\right]^\frac{1}{2} \right)^{\frac{1}{d+1}}. \\
\end{eqnarray*}
In other words we obtained that
\begin{equation}
\label{E[sup|X_{t,s}(x)|^gamma] <= ... (N=3)}
       E\left[ \sup_{|x|\leq m} |X_{t,s}^N(x)|^\gamma \right]
  \leq C_{m,\gamma}(s-t)^{\frac{N}{2}\gamma - \frac{1}{d+1}}, \quad
       0 \leq t < s \leq T.
\end{equation}

Next, recall the multiple integrals $X_{t,s}^{k,l}(x)$ defined by
(\ref{Xstkl}). For $0 \leq r_1 \leq r_2 \leq r_3 \leq T$, and $1\le
n_0\le N$, define
 \bea
 \label{YNrt}
  Y_{r_1,r_2,r_3,n_0}^N(x).
  \dfnn
 X_{r_2,r_3}^{N-n_0,N}(x) X_{r_1,r_2}^{1,N-n_0-1}(x),
 \eea
 An easy calculation shows that
\begin{eqnarray}
\label{Xrecurs}
      X_{r_1,r_3}^{1,N}(x)
\nonumber
& = & \int_{r_1}^{r_3} X_{r_1,u_N}^{1,N-1}(x) f_N(u_N,x) dB_{u_N} \\
\nonumber
& = & \int_{r_1}^{r_2} X_{r_1,u_N}^{1,N-1}(x) f_N(u_N,x) dB_{u_N}
      + \int_{r_2}^{r_3} X_{r_1,u_N}^{1,N-1}(x) f_N(u_N,x) dB_{u_N} \\
\label{X_{r_1,r_3}(x)}
& = & X_{r_1,r_2}^{1,N}(x) + \sum_{n_0=0}^{N-2} Y_{r_1,r_2,r_3,n_0}^N(x) + X_{r_2,r_3}^{1,N}(x)
\end{eqnarray}

To simplify the further discussion, we now assume without loss of
generality that $T=1$. Let $D_n$ be the set of all $t_i^n :=
2^{-n}i$, for some $0\leq i\leq 2^n$. The set $D=\bigcup_{n\geq 1}
D_n$ is then the set of all dyadic numbers in $[0,1]$. Now, let
$n\in\hN$ and $0 \leq t < s \leq 1$ be arbitrary dyadic numbers such
that $s-t\leq 2^{-n}$. Our aim is to estimate $X_{t,s}(x)$ uniformly
in $t,s\in D$. To this end, we notice that there exists some $k\geq
n$ such that $t,s$ belong to $D_k$. Moreover, denoting

\begin{figure}[h]
\begin{center}
\includegraphics[width=127.6mm, height=100mm, angle=0]{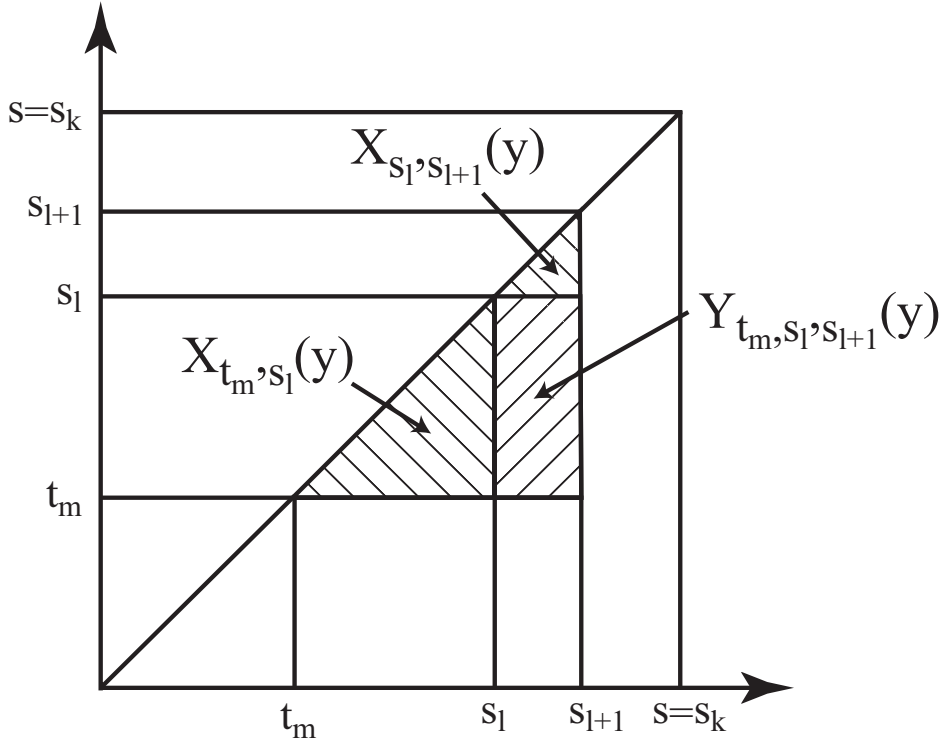}
\caption{$X_{s,t}=X^{1,2}_{s,t}$ and
$Y_{t,s_1,s_2}=Y^2_{t,s_1,s_2,n_0}$ for $N=2$ and $n_0=1$}
\label{Estimate_Iterated_Ito_Integral}
\end{center}
\end{figure}


\[
  s_j = \sup\{r\in D_j, \; r\leq s \}, \quad t_j = \sup\{r\in D_j, \; r\leq t \}, \quad n\leq j\leq k,
\]
we have $s_n \leq s_{n+1} \leq \ldots \leq s_k = s$, $t_n \leq
t_{n+1} \leq \ldots \leq t_k = t$. Therefore (see also Figure \ref{Estimate_Iterated_Ito_Integral})
\begin{eqnarray}
\label{X_{t_n,s} Eq. 1}
  X_{t_n,s}^{1,N}(x) & = & X_{t_n,s_n}^{1,N}(x)
                           + \sum_{j=n}^{k-1} \left\{ X_{s_j,s_{j+1}}^{1,N}(x) + \sum_{n_0=0}^{N-2} Y_{t_n,s_j,s_{j+1},n_0}^N(x)  \right\}, \\
\label{X_{t_n,s} Eq. 2}
  X_{t_n,s}^{1,N}(x) & = & X_{t,s}^{1,N}(x) + \sum_{j=n}^{k-1} \left\{ X_{t_j,t_{j+1}}^{1,N}(x) + \sum_{n_0=0}^{N-2} Y_{t_j,t_{j+1},s,n_0}^N(x) \right\}.
\end{eqnarray}
In order to estimate $Y_{t_n,s_j,s_{j+1},n_0}^N(x)$ and
$Y_{t_j,t_{j+1},s,n_0}^N(x)$, we notice that
\[
  s_{j+1} - s_j \in \{0,2^{-(j+1)}\}, \quad t_{j+1} - t_j \in \{0,2^{-(j+1)}\}, \quad
  \quad n\leq j\leq k,
\]
and thus
 \bea
 \label{Estimate for s_j - t_n}
 \left\{\begin{array}{lll}
  0 \leq s_j - t_n = (s_j - s_n) + (s_n - t_n) \leq 2\cdot
  2^{-n};\\
 0 \leq s - t_{j+1} = (s -t) + (t - t_{j+1}) \leq 2\cdot 2^{-n},
 \end{array}
 \right.
 \eea
for $n\leq j\leq k$. Moreover, recall that for every
$\alpha\in(0,\frac{1}{2})$, there exists a $\theta_\alpha \in
L^{\infty}_-(\Omega, {\cal F}, P)$ such that for all $l, l'\in \hN$
with $0\le l'-l\le N-2$,
\[
  P\left\{ \left| X_{t,t+h}^{l,l'} \right| \leq \theta_\alpha h^{(l'-l+1)\alpha},
  \; 0 \leq t < t+h \leq T, \; |x| \leq m \right\} = 1,
\]
thanks to the inductional hypothesis. Thus, we have by
(\ref{Estimate for s_j - t_n})
\begin{eqnarray}
         \sum_{n_0=0}^{N-2}|Y_{t_n,s_j,s_{j+1},n_0}^N(x)|
\nonumber
&\leq & \sum_{n_0=0}^{N-2} \theta_\alpha |s_{j+1} - s_j|^{(n_0+1)\alpha} \cdot \theta_\alpha |s_j - t_n|^{(N-n_0-1)\alpha} \\
\nonumber
& \leq & \sum_{n_0=0}^{N-2} \theta_\alpha^2 2^{-(n_0+1)(j+1)\alpha} 2^{-(N-n_0-1)(n-1)\alpha} \\
\label{Estimate |Y_j^{1,1}| + |Y_j^{1,2}|}
& \leq & 2^N \sum_{n_0=0}^{N-2} \theta_\alpha^2 2^{-((n_0+1)j + (N-n_0-1)n)\alpha}.
\end{eqnarray}
Similarly, we obtain
\begin{equation}
\label{Estimate |Y_j^{2,1}| + |Y_j^{2,2}|}
  \sum_{n_0=0}^{N-2} |Y_{t_j,t_{j+1},s,n_0}^N(x)| \leq 2^N \sum_{n_0=0}^{N-2}
\theta_\alpha^2 2^{-((n_0+1)j + (N-n_0-1)n)\alpha}.
\end{equation}
Next, we define
$  K_j \dfnn \sup\{ |X_{t_i^j,t_{i+1}^j}(x)|, \; |x|\leq m, \; 0
\leq i \leq 2^j-1 \}$, $n\leq j\leq k$.
Then, combining (\ref{X_{t_n,s} Eq. 1}), (\ref{X_{t_n,s} Eq. 2}),
(\ref{Estimate |Y_j^{2,1}| + |Y_j^{2,2}|}), and (\ref{Estimate
|Y_j^{1,1}| + |Y_j^{1,2}|}), we deduce
\begin{eqnarray*}
&  & |X_{t,s}^N(x)| \\
& =& \left| X_{t_n,s_n}^N(x) + \sum_{j=n}^{k-1} \left\{
X_{s_j,s_{j+1}}^N(x) - X_{t_j,t_{j+1}}^N(x)
         + \sum_{n_0=0}^{N-2} Y_{t_n,s_j,s_{j+1},n_0}^N(x) - \sum_{n_0=0}^{N-2} Y_{t_j,t_{j+1},s,n_0}^N(x) \right\} \right|, \\
&\leq& K_n + \sum_{j=n}^{k-1} \left( 2K_j + 2^{N+1} \sum_{n_0=0}^{N-2} \theta_\alpha^2 2^{-((n_0+1)j + (N-n_0-1)n)\alpha} \right) \\
&\le& \sum_{j=n}^\infty \left( 3K_j + 2^{N+1} \sum_{n_0=0}^{N-2}
\theta_\alpha^2 2^{-((n_0+1)j + (N-n_0-1)n)\alpha} \right).
\end{eqnarray*}

Now let $(N-1)\alpha < \beta < N\alpha < \frac{N}{2}$. Define
$$ \zeta^n_{\beta,m} \dfnn  \sup\left\{ \frac{|X_{t,s}^N(x)|}{|s-t|^\beta}, \; s,t\in D, \;
         2^{-(n+1)} < s-t \leq 2^{-n}, \; |x|\leq m  \right\},
$$
and $M_{\beta,m}\dfnn \sup_{n\geq 1}\zeta^n_{\beta,m}$. Then,
\begin{eqnarray*}
         M_{\beta,m}
& \leq & \sup_{n\geq 1} \left( 2^{(n+1)\beta} \sup\left\{
|X_{t,s}^N(x)|, \; s,t\in D, \;
         0 < s-t \leq 2^{-n}, \; |x|\leq m \right\} \right) \\
& \leq & \sup_{n\geq 1} \left( 2^{(n+1)\beta}
         \sum_{j=n}^\infty \left( 3K_j + 2^{N+1} \sum_{n_0=0}^{N-2} \theta_\alpha^2 2^{-((n_0+1)j + (N-n_0-1)n)\alpha} \right) \right) \\
& \leq & \sup_{n\geq 1} 2^{\beta+N+2} \sum_{j=n}^\infty \left( 2^{n\beta} K_j +
         \sum_{n_0=0}^{N-2} \theta_\alpha^2 2^{-(n_0+1)\alpha j + (\beta-(N-n_0-1)\alpha)n}
         \right) \\
& \leq & \sup_{n\geq 1} 2^{\beta+N+2} \sum_{j=n}^\infty
\left(2^{n\beta} K_j + \sum_{n_0=0}^{N-2} \theta_\alpha^2
2^{-(n_0+1)\alpha j + (\beta-(N-1)\alpha)n+nn_0\a} \right) \\
& \leq & 2^{\beta+N+2} \sum_{j=0}^\infty \left( 2^{j\beta} K_j +
         (N-1)\theta_\alpha^2 2^{(\beta-N\alpha)j} \right) \\
& =    & 2^{\beta+N+2} \left( \sum_{j=0}^\infty 2^{j\beta} K_j + \frac{(N-1)\theta_\alpha^2}{1-2^{-(N\alpha-\beta)}} \right).
\end{eqnarray*}

Noting that by (\ref{E[sup|X_{t,s}(x)|^gamma] <= ... (N=3)})
\begin{equation}
\label{Estimate for E[K_j^gamma]}
       E[K_j^\gamma]
  \leq \sum_{i=0}^{2^j-1} E\left[ \sup_{|x|\leq m} |X_{t_i^j,t_{i+1}^j}(x)|^\gamma \right]
  \leq 2^jC_{m,\gamma}2^{-j(\frac{N}{2}\gamma-1)}
  =    C_{m,\gamma}2^{-j(\frac{N}{2}\gamma-2)},
\end{equation}
we thus obtain for $\gamma>\frac{4}{N - 2\beta}$ that
\begin{eqnarray*}
         \left( E\left[ M_{\beta,m}^\gamma \right] \right)^{\frac{1}{\gamma}}
& \leq & 2^{\beta+N+2} \left( \sum_{j=0}^\infty 2^{j\beta} \left( E\left[ K_j^\gamma \right] \right)^{\frac{1}{\gamma}} +
         \frac{N-1}{1-2^{-(N\alpha-\beta)}} \left( E\left[ \theta_\alpha^{2\gamma} \right] \right)^{\frac{1}{\gamma}}
         \right) \\
& \leq & 2^{\beta+N+2} \left( \sum_{j=0}^\infty 2^{j\beta} C_{m,\gamma}2^{-j
         \left( \frac{N}{2} - \frac{2}{\gamma} \right) } +
         \frac{N-1}{1-2^{-(N\alpha-\beta)}} \left( E\left[ \theta_\alpha^{2\gamma} \right] \right)^{\frac{1}{\gamma}}
         \right) \\
& =    & 2^{\beta+N+2} C_{m,\gamma} \left( \frac{1}{1-2^{-\left( \frac{N}{2} - \frac{2}{\gamma}-\beta \right)}} +
         \frac{N-1}{1-2^{-(N\alpha-\beta)}} \left( E\left[ \theta_\alpha^{2\gamma} \right] \right)^{\frac{1}{\gamma}}
         \right) < \infty.
\end{eqnarray*}
Consequently, $M_{\beta,m}\in L^{\infty}_-(\Omega,{\cal F}, P)$, for
all $\beta \in (0,N\alpha)$. But since $\alpha\in(0,\frac{1}{2})$ is
arbitrary, we can extend the result to $\beta\in(0,N/2)$. Finally,
note that the definitions of $M_{\beta,m}$ and
$\{\zeta^n_{\beta,m}\}$, as well as the continuity of the mapping
$(t,s) \mapsto X_{t,s}(x)$, imply that
$$
\h\zeta_{\beta,m}\dfnn \sup\left\{
\frac{|X_{t,s}^N(x)|}{|s-t|^\beta}, \;   0< s-t \leq 2^{-1},
\;|x|\leq m  \right\}\in L^{\infty}_-(\Omega,{\cal F}, P),
 $$
for all $\beta\in(0,N/2)$ and $m\in\hN$.  The proposition then
follows from the recursive relation (\ref{Xrecurs}).\qed


The following corollary can be easily obtained by adapting the proof
of Proposition \ref{Estimation iterated Ito integral} in an obvious
manner.
\begin{cor}
\label{Cor. to Estimation iterated Ito integral} The statement of
Proposition \ref{Estimation iterated Ito integral} remains valid if,
for $1\leq i\leq N$, $dB_{s_i}$ is replaced by $dA^i_{s_i}$ in
(\ref{Xstj}), where $A^i$ is either the Brownian motion $B$ or
$A^i_s=s$, $s\in[0,T]$. Moreover, (\ref{zeta_{beta,m} in L^infty-})
holds whenever $\beta \in \big(0,l_1+\frac{l_2}{2}\big)$,
$m\in\dbN$, where $l_1$ is the number of $i$, $1\leq i\leq N$, for
which $A^i_s=s$, $s\in[0,T]$, and $l_2=N-l_1$.
\end{cor}

\section{Forward Taylor expansion}
\setcounter{equation}{0}

In this and next section we shall provide a complete proof of
Theorem \ref{StoTaylor}. The results in each section, however, can
also be applied independently, and therefore are of interest in
their own right. We should note that, unlike the usual Taylor
expansion, in the stochastic case the direction of the time
increment makes significant difference in the argument, due to the
``progressive measurability" of the random fields. We will thus
separate the two cases: in this section we study the {\it forward}
expansion, and leave the {\it backward} case to the next section.

\subsection{Forward temporal expansion}

We begin with the forward ``temporal" expansion, that is, only the
time variable has the increment. Let us first introduce the
following extra notations. For each $T>0$ we define the
2-dimensional simplex
 \bea
 \label{simplex}
  \triangle_{[0,T]} \dfnn \{(t,s)\in[0,T]\times\hR^+:0 \leq t < s \leq
  T\}.
 \eea
For any real-valued measurable functional $\theta$ defined on
$\Omega \times \triangle_{[0,T]} \times \dbR^{d+l}$ and $m\in\hN$,
we denote
 \bea
 \label{thetahat}
     \widehat\theta_m
     \dfnn \sup\{ |\theta(t,s,x,z)|, \; (t,s)\in\triangle_{[0,T]}, ~x\in\hR^d, \; z\in\hR
     \mbox{ with } |x|,|z| \leq m \}.
 \eea
Furthermore, by a slight abuse of notations, in what follows we
shall denote $R_\a$ to be any measurable functional $\theta_\a$,
indexed by $\alpha\in\left(\frac{1}{3},\frac{1}{2}\right)$,  such
that $(\widehat\theta_\a)_m\in L^{\infty}_-(\Omega,\cF,P)$ for all
$m\in\hN$, and again, it may vary from line to line.

Our main result of this section is the following  {\it Stochastic
forward temporal Taylor expansion}.
\begin{prop}
\label{zeta(t+h,x) - zeta(t,x)} Assume that $\zeta\in
C^{0,(3)}(\bF^B, [0,T]\times\hR^d)$ satisfies (T1)-(T3). Then, for
all $m\in\dbN$ and $\alpha\in\left(\frac{1}{3},\frac{1}{2}\right)$
there exists a subset $\widetilde\Omega \subset \Omega$ such that
$P(\widetilde\O)=1$, and that on $\widetilde\Omega$, for all $0
\leq t < t+h \leq T$, $x\in\hR^d$, the following expansion holds
 \bea
 \label{ftempexp}
 \zeta(t+h,x) - \zeta(t,x)= ah + b(B_{t+h} - B_t)+ \frac{c}{2} (B_{t+h} - B_t)^2
      + h^{1+\alpha}R_\alpha(t,t+h,x),
 \eea
where
 \bea
 \label{ftempexp1}
 a &=& F(x,\zeta_2(t,x))- \frac{1}{2} (D_z
 G)(x,\zeta_2(t,x))G_2(x,\zeta_3(t,x)),\nonumber\\
 b &=&
G(x,\zeta_2(t,x)),\\
 c&=&(D_zG)(x,\zeta_2(t,x))G_2(x,\zeta_3(t,x)). \nonumber
 \eea
\end{prop}

{\it Proof.} Let $h>0$ be such that $0\leq t \leq t+h \leq T$. For
any $x\in\hR^d$, we write
 \bea
 \label{dhzeta+}
  \zeta(t+h,x) - \zeta(t,x)  =  \int_t^{t+h} F(x,\zeta_2(s,x)) ds + \int_t^{t+h} G(x,\zeta_2(s,x)) dB_s
                            \dfnn & I^1+ I^2,
 \eea
where $I^i=I^i(t,h,x)$, $i=1,2$, are the two integrals. We shall
study their expansions separately.

We begin by $I^1$. The argument is very similar to that of
\cite{Buckdahn_Ma_02}, we provide a sketch for completeness.
Let
 \beaa \label{Def. tilde(H)}
  H^1(x,z_2,z_3) \dfnn (D_zF)(x,z_2)F_2(x,z_3) + \frac{1}{2}\tr[G_2(G_2)^*(x,z_3)D_z^2F(x,z_2)]
 \eeaa
for $(x,z_2,z_3)\in\hR^d\times\hR^{d_2}\times\hR^{d_3}$. Then,
applying  It\^o's formula, we have
 \bea
 \label{DeltaF}
 F(x,\zeta_2(s,x)) - F(x,\zeta_2(t,x)) & = &
 \int_t^sH^1(x,\zeta_2(r,x),\zeta_3(r,x)) dr \nonumber\\
 &   & + \int_t^s \lan(D_zF)(x,\zeta_2(r,x)),G_2(x,\zeta_3(r,x))\ran dB_r,
 \eea
for $0 \leq t \leq s \leq T$ and $x\in\hR^d$, $P$-a.s. We now show
that there is a universal subset $\Omega'\subseteq \Omega$ with
$P(\O')=1$, on which (\ref{DeltaF}) holds for all $0\leq t\leq s\leq
T$. But for this it suffices to prove the random fields
$$(s,t,x) \mapsto \int_t^s
H^1(x,\zeta_2(r,x),\zeta_3(r,x)) dr, \q \int_t^s
\lan(D_zF)(x,\zeta_2(r,x)),G_2(x,\zeta_3(r,x)) \ran dB_r
$$
have continuous versions. To see this, we try to make use of the
Kolmogorov continuity criterion.
For notational simplicity let us denote
$$ \G(t,x)\dfnn (D_zF)(x,\zeta_2(t,x))G_2(x,\zeta_3(t,x)), \qq(t,x)\in [0,T]\times
\hR^d.
$$
Then, we can deduce from (T3) and the Burkholder-Davis-Gundy
inequality that, for each $k>2$,
\begin{eqnarray*}
&      & E\bigg| \int_t^s \G(r,x) dB_r  -\int_{t'}^{s'}
\G(r,x')dB_r \bigg|^k \\
& \leq & C_k\left\{E\Big| \int_t^{t'} |(\G(r,x)|^2 dr
         \Big|^\frac{k}{2}+E\Big| \int_{s'}^s
|\G(r,x)|^2 dr\Big|^\frac{k}{2}
+E\bigg( \int_0^T |\G(r,x)-\G(r,x') |^2 dr \bigg)^\frac{k}{2} \right\} \\
& \leq & C_k\left( |t'-t|^\frac{k}{2} + |s'-s|^\frac{k}{2} \right)
E\bigg\{ \sup_{r\in[0,T], \atop |x|\leq m}
         |\G(r,x)|^k \bigg\}+C_k|x-x'|^k E\bigg\{ \sup_{r\in[0,T], \atop |x|\leq m} |D_x
         \G(r,x)|^k \bigg\} \\
& \leq & C_k\left( |t'-t|^\frac{k}{2} + |s'-s|^\frac{k}{2} +
|x-x'|^k \right)
\end{eqnarray*}
for $0\leq s \leq t \leq T$, $0\leq s' \leq t' \leq T$,
$|x|,|x'|\leq m $, and $k\in\hN$. Hence, the Kolmogorov continuity
criterion renders that the random field $\Big\{ \int_t^s
\G(r,x)dB_r, \; 0\leq t\leq s\leq T, \; |x|\leq m \Big\}$ possesses
a version that is continuous in $(t,s,x)$.
A similar estimate allows to prove that also the first integral in
(\ref{DeltaF}) admits a version continuous in $(t,s,x)$ as well.
Hence, we conclude that on some $\Omega'\subset\Omega$ of full
probability measure the relation (\ref{DeltaF}) holds for all $0\leq
t\leq s\leq T$ and $|x|,|x'|\leq m$.

Consequently, writing the integral $I^1$ as
\begin{eqnarray*}
      I^1(t,h,x)
& = & F(x,\zeta_2(t,x))h + \int_t^{t+h} \int_t^s {H}^1(x,\zeta_2(r,x),\zeta_3(r,x)) dr ds \\
&   & + \int_t^{t+h} \int_t^s
(D_zF)(x,\zeta_2(r,x))G_2(x,\zeta_3(r,x)) dB_r ds,
\end{eqnarray*}
and using the conclusion above we see that as a function of
$(t,h,x)$, $I^1$ is jointly continuous for all $0\leq t\leq t+h\leq
T$ and $x\in\hR^d$, over the universal set $\O'$. It then follows
from Corollary \ref{Cor. to Estimation iterated Ito integral} that
there exists an $\O''\subseteq \O'$ with $P(\O'')=1$ such that for
all $0\leq t\leq t+h\leq T$ and $x\in\hR^d$, it holds that
 \begin{equation}
 \label{Estimate Delta^1}
  I^1(t,h,x) = F(x,\zeta_2(t,x)h + h^{1 + \alpha}R_\alpha(t,t+h,x).
 \end{equation}
In what follows we will not distinguish $\O''$ from $\O'$.

We now turn our attention to $I^2$. Again, we begin by denoting,
for $(x,z_2,z_3) \in \hR^d\times\hR^{d_2}\times\hR^{d_3}$,
\begin{equation}
\label{Def. H}
  H^2(x,z_2,z_3) \dfnn (D_zG)(x,z_2)F_2(x,z_3) + \frac{1}{2}\tr[G_2(G_2)^*(x,z_3)
D_z^2G(x,z_2)].
\end{equation}
Then, for every $x\in\hR^d$, we again apply It\^o's formula to
get, for $0 \leq t \leq s \leq T$, and $P$-a.s.,
\begin{eqnarray}
\label{dGzeta}
 G(x,\zeta_2(s,x)) - G(x,\zeta_2(t,x))
& = & \int_t^s H^2(x,\zeta_2(r,x),\zeta_3(r,x)) dr \\
&   & + \int_t^s (D_zG)(x,\zeta_2(r,x))G_2(x,\zeta_3(r,x))
dB_r.\nonumber
\end{eqnarray}
Using the similar arguments as before we can find another universal
subset, still denoted by $\Omega'\subseteq \Omega$ with $P(\O')=1$,
on which (\ref{dGzeta}) holds for all $0\leq t\leq s\leq T$.

Next, using (\ref{dGzeta}) it is easy to see that $I^2$ can be
written as
\begin{eqnarray}
 \label{Eq. for Delta^2}
      I^2(t,h,x)
\nonumber
& =  & G(x,\zeta_2(t,x))(B_{t+h} - B_t) + \int_t^{t+h} \int_t^s H^2(x,\zeta_2(r,x),\zeta_3(r,x)) dr dB_s \\
\nonumber
&    & + \int_t^{t+h} \int_t^s (D_zG)(x,\zeta_2(r,x))G_2(x,\zeta_3(r,x)) dB_r dB_s \\
& \dfnn & G(x,\zeta_2(t,x))(B_{t+h} - B_t) + I^{2,1}(t,h,x) +
I^{2,2}(t,h,x),
\end{eqnarray}
and we can claim as before that (\ref{Eq. for Delta^2}) holds for
all $0\leq t<t+h\leq T$, $|x|\leq m$, over a
universal subset of full probability measure, again denoted by
$\Omega'$.

 We now analyze $I^{2,1}$ and $I^{2,2}$ separately.
Using integration by parts we see that
 \beaa
       I^{2,1}
 = (B_{t+h} - B_t)\neg \int_t^{t+h}\neg H^2(x,\zeta_2(r,x),\zeta_3(r,x)) dr
 -\neg  \int_t^{t+h}\neg (B_r\neg - \neg B_t) H(x,\zeta_2(r,x),\zeta_3(r,x)) dr.
 \eeaa
Then, following the argument developed in the previous part, we can
show that the equality holds for all $0\leq t < t+h \leq T$,
$x\in\hR^d$, over an $\O'\subseteq \O$ with $P(\O')=1$. Note that
\begin{equation}
\label{Hoelder continuity of BM}
  \sup_{0 \leq t < t+h \leq T}\{ h^{-\alpha}|B_{t+h} - B_t| \} \in L^{\infty}_-
  (\O, {\cal F}, P),
\end{equation}
it follows that over $\O'$,
\begin{equation}
\label{Estimate Delta^2,1}
  I^{2,1}(t,h,x) = h^{1+\alpha}R_\alpha(t,t+h,x).
\end{equation}

The estimate for $I^{2,2}$ is slightly more involved. For notational
simplicity let us define
\begin{eqnarray*}
            \widehat{K}(x,z_2,z_3)
  & \dfnn & (D_zG)(x,z_2)G_2(x,z_3), \quad (x,z_2,z_3) \in \hR^d\times
\hR^{d_2}\times\hR^{d_3}, \\
         \widehat{F}(x,z_3,z_4)
  & \dfnn & (F_2(x,z_3),F_3(x,z_4)), \quad (x,z_3,z_4) \in
\hR^d\times\hR^{d_3}\times\hR^{d_4}, \\
         \widehat{G}(x,z_3,z_4)
  & \dfnn& (G_2(x,z_3),G_3(x,z_4)),  \\
         \widehat{H}(x,\widehat{z}_2,\widehat{z}_3)
  & \dfnn & D_{\widehat{z}_2}\widehat{K}(x,\widehat{z}_2)\widehat{F}(x,
\widehat{z}_3)
         + \frac{1}{2}\tr[\widehat{G}\widehat{G}^*(x,\widehat{z}_3)D_{
\widehat{z}_2}^2
         \widehat{K}(x,\widehat{z}_2)],
\end{eqnarray*}
where $ \widehat{z}_2 = (z_2,z_3)$ and $\widehat{z}_3 =
(z_3,z_4)$. Moreover, we denote $  \widehat\zeta_i(s,x) =
(\zeta_i(s,x),\zeta_{i+1}(s,x))$, $i=2,3$.

 Then, applying  It\^o's formula we have,
\begin{eqnarray*}
      \widehat{K}(x,\widehat\zeta_2(s,x)) - \widehat{K}(x,\widehat\zeta_2(t,x))
& = & \int_t^s \widehat{H}(x,\widehat\zeta_2(r,x),\widehat\zeta_3(r,x)) dr \\
&   & + \int_t^s
D_{\widehat{z}_2}\widehat{K}(x,\widehat\zeta_2(r,x))\widehat{G}(x,\widehat\zeta_3(r,x))
dB_r.
\end{eqnarray*}
Again, we assume that the equality holds
for all $0 \leq s \leq t \leq T$ and $x\in\hR^d$, on $\O'$.
Therefore
\begin{eqnarray}
\label{Eq. Delta^2,2}
I^{2,2}(t,h,x)
& = & \int_t^{t+h}\int_t^s \widehat{K}(x,\widehat\zeta_2(r,x))dB_r dB_s \nonumber \\
& = & \frac{1}{2}(D_zG)(x,\zeta_2(t,x))G_2(x,\zeta_3(t,x))(B_{t+h} - B_t)^{\diamond 2} \\
\nonumber
&   & + \int_t^{t+h} \int_t^s \int_t^r \widehat{H}(x,\widehat\zeta_2(v,x),\widehat\zeta_3(v,x))
dv dB_r dB_s \\
 &   & + \int_t^{t+h} \int_t^s \int_t^r
D_{\widehat{z}_2}\widehat{K}(x,\widehat\zeta_2(v,x))\widehat{G}(x,\widehat\zeta_3(v,x))
dB_v dB_r dB_s. \nonumber
\end{eqnarray}
Moreover,  Proposition \ref{Estimation iterated Ito integral}
implies that
\begin{eqnarray*}
h^{-(1+\alpha)} \int_t^{t+h} \int_t^s \int_t^r
      D_{\widehat{z}_2}\widehat{K}(x,\widehat\zeta_2(v,x)) \widehat{G}(x,\widehat\zeta_3(v,x))
      dB_v dB_r dB_s = R_\alpha(t,t+h,x);
\end{eqnarray*}
and  Corollary \ref{Cor. to Estimation iterated Ito integral}
implies that
\[
   h^{-(1+2\alpha)} \int_t^{t+h} \int_t^s \int_t^r
   \widehat{H}(x,\widehat\zeta_2(v,x),\widehat\zeta_3(v,x)) dv dB_r dB_s
 = R_\alpha(t,t+h,x).
\]
Consequently, we see that, over $\O'$, (\ref{Eq. Delta^2,2}) becomes
\begin{equation}
\label{Estimate Delta^2,2}
    I^{2,2}(t,h,x)
  = \frac{1}{2}(D_zG)(x,\zeta_2(t,x))G_2(x,\zeta_3(t,x))(B_{t+h} - B_t)^{\diamond 2}
    + h^{1+\alpha}R_\alpha(t,t+h,x).
\end{equation}

Finally, plugging (\ref{Estimate Delta^1}), (\ref{Eq. for Delta^2}),
(\ref{Estimate Delta^2,1}), and (\ref{Estimate Delta^2,2}) into
(\ref{dhzeta+}),
 we obtain
 (\ref{ftempexp}) and (\ref{ftempexp1}), with a universal exceptional null set,
 proving the proposition. \qed

\subsection{Forward temporal-spatial Taylor expansion}

Based on the forward temporal Taylor expansion  Proposition
\ref{zeta(t+h,x) - zeta(t,x)}, we now add the spatial increment. Our
main result of this section is the following proposition.

\begin{prop}
\label{zeta(t+h,x+k) - zeta(t,x)} Assume that $\zeta\in
C^{0,(3)}(\bF^B, [0,T],\hR^d)$ satisfying (T1)-(T3). Then, for all
$m\in\dbN$ and $\alpha\in\left(\frac{1}{3},\frac{1}{2}\right)$ there
exists some subset $\widetilde\Omega \subset \Omega$ such that
$P(\widetilde\O)=1$, and that on $\widetilde\Omega$, for all $0 \leq
t < t+h \leq T$ and $x,k\in\hR^d$,
 \bea
 \label{dhxzeta+}
 &&\zeta(t+h,x+k) - \zeta(t,x)= ah+b(B_{t+h} - B_t) + \frac{c}{2} (B_{t+h} - B_t)^2
      + \lan p,k \ran+ \frac{1}{2}\langle Xk,k \rangle \nonumber\\
 && \qq\qq\qq\qq\qq\q+\lan q, k\ran (B_{t+h} - B_t) + (h
 +|k|^2)^{3\alpha}R_\alpha(t,t+h,x,x+k),
 \eea
where
 \bea
 \label{dhxzeta+1}
 a&=&  F(x,\zeta_2(t,x)) - \frac{1}{2} (D_z
G)(x,\zeta_2(t,x))G_2(x,\zeta_3(t,x));\nonumber\\
 b&=& G(x,\zeta_2(t,x));\q
 c= (D_z G)(x,\zeta_2(t,x))G_2(x,\zeta_3(t,x));\\
 p&=& (D_x\zeta)(t,x); \q
 X= D_x^2\zeta(t,x) \nonumber;\\
 q&=& (D_x G)(x,\zeta_2(t,x)) + (D_z G)(x,\zeta_2(t,x))D_x\zeta_2(t,x)\nonumber .
 \eea

\end{prop}


{\it Proof.}
First let us write
 $$\zeta(t+h,x+k) - \zeta(t,x)=[\zeta(t+h,x+k) -
 \zeta(t+h,x)]+[\zeta(t+h,x) - \zeta(t,x)],
 $$
where the second $[\cds]$ above is the forward temporal expansion
studied in the previous subsection. In light of Proposition
\ref{zeta(t+h,x) - zeta(t,x)}, we need only prove the following {\it
stochastic spatial Taylor expansion}:
 \bea
 \label{stospexp}
  \zeta(t+h,x+k) - \zeta(t+h,x) & =& D_x\zeta(t,x)k + \frac{1}{2}\langle D_x^2\zeta(t,x)k,k \rangle
 + \{(D_xG)(x,\zeta_2(t,x)) \nonumber\\
 & &+ (D_zG)(x,\zeta_2(t,x))D_x\zeta_2(t,x)\}k(B_{t+h} - B_t)\\
 && + (h + |k|^2)^{3\alpha} R^1_\alpha(t,t+h,x,x+k). \nonumber
 \eea

To this end, we first apply the standard Taylor expansion and use
the assumption (T3) to get
\begin{eqnarray}
      \zeta(t+h,x+k) - \zeta(t+h,x)
\nonumber
& = & D_x\zeta(t+h,x)k + \frac{1}{2}\langle D_x^2\zeta(t+h,x)k,k \rangle \\
\label{Eq. zeta(t+h,x+k) - zeta(t+h,x)} &   & +
|k|^3R^1_\a(t,t+h,x,x+k),
\end{eqnarray}
for all $0 \leq t < t+h \leq T$, $x,k\in\hR^d$, $P$-a.s. Next,
differentiating the equation for $\zeta=\zeta_1$ in (\ref{zetai}) we
have
\begin{eqnarray*}
 D_x\zeta(t,x)
& = & D_x\zeta_0(x) + \int_0^t \left\{ (D_xF)(x,\zeta_2(s,x))
      + (D_zF)(x,\zeta_2(s,x))D_x\zeta_2(s,x) \right\} ds \\
&   & + \int_0^t \left\{ (D_xG)(x,\zeta_2(s,x))
      + (D_zG)(x,\zeta_2(s,x))D_x\zeta_2(s,x)) \right\} dB_s,
\end{eqnarray*}
for all $t\in[0,T]$, $x\in\hR^d$, $P$-a.s. Now,
applying Proposition \ref{zeta(t+h,x) - zeta(t,x)} (to $D_x\zeta$)
one can check that
\begin{eqnarray*}
&&D_x\zeta(t+h,x)\\
& = & D_x\zeta(t,x) + \{(D_xF)(x,\zeta_2(t,x)) + (D_zF)(x,\zeta_2(t,x))D_x\zeta_2(t,x)\}h \\
&   & -\frac12\{D_x[(D_zG)(x,
\zeta_2(t,x))]G_2(x,\zeta_3)+D_zG(x,\zeta_2(t,x))D_x[G_2(x,\zeta_3)]\}h\\
&   & + \{(D_xG)(x,\zeta_2(t,x)) + (D_zG)(x,\zeta_2(t,x))D_x\zeta_2(t,x)\}(B_{t+h} - B_t) \\
&   & + \frac{1}{2} \{D_x[(D_zG)(x,\zeta_2(t,x))]G_2(x,\zeta_3)
      + D_zG(x,\zeta_2(t,x))D_x[G_2(x,\zeta_3)]\}(B_{t+h} - B_t)^{\diamond 2} \\
&   & + h^{1+\alpha}R_\alpha(t,t+h,x),
\end{eqnarray*}
for all $0 \leq t < t+h \leq T$, $x\in\hR^d$, $P$-a.s.
Consequently, it follows from  (\ref{Hoelder continuity of BM}) that
\begin{eqnarray}
\label{Eq. D_x zeta(t+h,x)}
D_x\zeta(t+h,x)
& = & D_x\zeta(t,x) + \{(D_xG)(x,\zeta_2(t,x)) + (D_zG)(x,\zeta_2(t,x))D_x\zeta_2(t,x)\}(B_{t+h} - B_t)
\nonumber\\
&   &  + h^{2\alpha}R_\alpha(t,t+h,x),
\end{eqnarray}
for all $0 \leq t < t+h \leq T$, $x\in\hR^d$, $P$-a.s. Similarly,
one shows that
\begin{equation}
\label{D_x^2 zeta(t+h,x)}
  D_x^2\zeta(t+h,x) = D^2_x\zeta(t,x) + h^\alpha R_\alpha(t,t+h,x), \quad
  0 \leq t < t+h \leq T, \; x\in\hR^d.
\end{equation}
Combining (\ref{Eq. zeta(t+h,x+k) - zeta(t+h,x)}), (\ref{Eq. D_x
zeta(t+h,x)}), and (\ref{D_x^2 zeta(t+h,x)}), we obtain that
\begin{eqnarray*}
&   & \zeta(t+h,x+k) - \zeta(t+h,x) \\
& = & D_x\zeta(t+h,x)k + \frac{1}{2}\langle D_x^2\zeta(t+h,x)k,k
\rangle
      + |k|^3R^1_\a(t,t+h,x,x+k) \\
& = & D_x\zeta(t,x)k + \frac{1}{2}\langle D_x^2\zeta(t,x)k,k \rangle \\
&   & + \{(D_xG)(x,\zeta_2(t,x)) + (D_zG)(x,\zeta_2(t,x))D_x\zeta_2(t,x)\}k(B_{t+h} - B_t) \\
&   & + |k|^3R^1_\a(t,t+h,x,x+k)+ h^{2\alpha}R_\alpha(t,t+h,x)k
      + h^\alpha\langle R_\alpha(t,t+h,x)k,k \rangle,
\end{eqnarray*}
for $0 \leq t < t+h \leq T$, $x,k\in\hR^d$, $P$-a.s. Finally, noting
that
\[
       |k|^3 + h^{2\alpha}|k| + h^\alpha|k|^2
  \leq C_m(h^{3\alpha} + |k|^3)
  \leq C_m(h + |k|^2)^{3\alpha}
\]
for all $h\in[0,T]$, $|k| \leq m$, we see that (\ref{stospexp})
holds, hence the proposition follows. \qed

\section{Backward Taylor expansion}
\setcounter{equation}{0}

In this section we treat the backward Taylor expansion, that is,
when the temporal increments are negative. As a general belief such
an expansion would be more difficult than the forward one, due to
the obvious ``adaptedness" issue. But we shall see, with our
``pathwise" approach, such difficulty is eliminated. We nevertheless
would like to separate its proof from the forward case because of
the slight difference in the arguments.
We again take two steps: first the backward
temporal expansion, and then the mixed time-space expansion.

\subsection{Backward temporal expansion}

We have the following analogy of Proposition \ref{zeta(t+h,x) -
zeta(t,x)}.
\begin{prop}
\label{zeta(t-h,x) - zeta(t,x)} Assume that $\zeta\in
C^{0,(3)}(\bF^B,[0,T],\hR^d)$ satisfies (T1)-(T3). Then, for all
$m\in\hN$ and $\alpha\in\left(\frac{1}{3},\frac{1}{2}\right)$, there
exists some subset $\widetilde\Omega \subset \Omega$ with
$P(\widetilde\O)=1$, such that on $\widetilde\Omega$, for all $0
\leq t-h < t \leq T$, $x\in\hR^d$, the following expansion holds
 \bea
 \label{dhzeta-}
 \zeta(t-h,x) - \zeta(t,x)= ah+ b(B_{t-h} - B_t) + \frac{c}{2} (B_{t-h} - B_t)^2
      + h^{1+\alpha}R_\alpha(t-h,t,x).
 \eea
where
 \bea
 \label{dhzeta-1}
 a&=& -\left\{ F(x,\zeta_2(t,x)) - \frac{1}{2} (D_z G)(x,\zeta_2(t,x))G_2(x,\zeta_3(t,x)) \right\};
 \nonumber\\
 b&=& G(x,\zeta_2(t,x));\\
 c&=& (D_z G)(x,\zeta_2(t,x))G_2(x,\zeta_3(t,x).\nonumber
 \eea
\end{prop}

{\it Proof.} Applying the forward Taylor expansion Proposition
\ref{zeta(t+h,x) - zeta(t,x)}, we have
\begin{eqnarray}
\nonumber \zeta(t,x) - \zeta(t-h,x)  \nonumber & = & \left\{
F(x,\zeta_2(t-h,x))
      - \frac{1}{2} (D_z G)(x,\zeta_2(t-h,x))G_2(x,\zeta_3(t-h,x))\right\}h \\
\nonumber
&   & + G(x,\zeta_2(t-h,x))(B_t - B_{t-h}) \\
\nonumber
&   & + \frac{1}{2} (D_z G)(x,\zeta_2(t-h,x))G_2(x,\zeta_3(t-h,x))(B_t - B_{t-h})^2 \\
\label{zeta(t,x) - zeta(t-h,x)} &   & +
h^{1+\alpha}R_\alpha(t-h,t,x),
\end{eqnarray}
for all $0 \leq t-h < t \leq T$, $x\in\hR^d$, on a full probability
set $\O'$. Our main task is to replace the temporal variable $t-h$
by $t$. To do this, we first apply It\^o's formula to
$G(x,\zeta_2(t,x))$ to obtain
 \bea
 \label{Gzeta2}
  G(x,\zeta_2(t,x)) = G(x,\zeta_2(0,x)) + \int_0^t H(x,\widehat\zeta_2(s,x)) ds
                      + \int_0^t L(x,\widehat\zeta_2(s,x)) dB_s,
 \eea
where $\widehat\zeta_2(t,x)=(\zeta_2(t,x),\zeta_3(t,x))$, and
\begin{equation}
  \label{Hz23}
  \left\{\begin{array}{lll}
  H(x,(z_2,z_3)) \dfnn  (D_zG)(x,z_2)F_2(x,z_3)
                               + \frac{1}{2}\tr[G_2(G_2)^*(x,z_3)(D_z^2G)(x,z_2)], \\
  L(x,(z_2,z_3)) \dfnn  (D_zG)(x,z_2)G_2(x,z_3).
  \end{array}\right.
\end{equation}
Next, we denote $\h z=(z_2,z_3)$, and define
$\widehat\zeta_3(t,x)\dfnn (\zeta_3(t,x),\zeta_4(t,x))$, $\widehat
G(t,x)\dfnn G(x,\zeta_2(t,x))$ and
$$ \widehat{G}^3(x,(z_3,z_4))=(G_2(x,z_3),G_3(x,z_4)).
$$
Applying Proposition \ref{zeta(t+h,x) - zeta(t,x)} to $ G(\cd,\cd)$
and using (\ref{Gzeta2}), we deduce that
\begin{eqnarray*}
  \h G(t,x) - \h G(t-h,x)& = & \left\{ H(x,\widehat\zeta_2(t-h,x))
      - \frac{1}{2} D_{\h z}L(x,\widehat\zeta_2(t-h,x))\widehat{G}^3(x,\widehat\zeta_3(t-h,x)) \right\}h \\
&   & + L(x,\widehat\zeta_2(t-h,x))(B_t - B_{t-h}) \\
&   & + \frac{1}{2} D_{\h z}L(x,\widehat\zeta_2(t-h,x))\widehat{G}^3(x,\widehat\zeta_3(t-h,x))(B_t - B_{t-h})^2 \\
&   & +  h^{1+\alpha}R_\alpha(t-h,t,x),
\end{eqnarray*}
for all $0 \leq t-h < t \leq T, \; x \in \hR^d$, which also holds on
the set $\O'$. Consequently, we obtain that
 \begin{eqnarray}
 \label{diffG}
 \h G(t,x) -\h G(t-h,x)  & = &
 L(x,\widehat\zeta_2(t-h,x))(B_t - B_{t-h}) +
 h^{2\alpha}R_\alpha(t-h,t,x) \nonumber\\
 & = & (D_zG)(x,\zeta_2(t-h,x))G_2(x,\zeta_3(t-h,x))(B_t - B_{t-h}) \\
 &   & +
 h^{2\alpha}R_\alpha(t-h,t,x). \nonumber
 \end{eqnarray}
In particular, by virtue of the H\"older continuity of the Brownian
motion (\ref{Hoelder continuity of BM}) and the assumption (T3) we
see from (\ref{diffG}) that
 \bea
 \label{G(x,zeta^2(t,x)) - G(x,zeta^2(t-h,x))}
  G(x,\zeta_2(t,x)) - G(x,\zeta_2(t-h,x)) = h^\alpha R_\alpha(t-h,t,x).
 \eea
Similarly, we can also derive the following (recall (\ref{Hz23})):
 \bea
 \label{F(x,zeta^2(t,x)) - F(x,zeta^2(t-h,x))}
 \left\{\begin{array}{lll}
  F(x,\zeta_2(t,x)) - F(x,\zeta_2(t-h,x)) = h^\alpha R_\alpha(t-h,t,x);\\
  L(x,\widehat\zeta_2(t,x))-L(x,\widehat\zeta_2(t-h,x))
 =  h^\alpha R_\alpha(t-h,t,x),
 \end{array}
 \right.
 \eea
for all $0\le t-h\leq t\leq T$, $|x|\leq m$, on the set $\O'$. Now
combining (\ref{Hz23})---(\ref{F(x,zeta^2(t,x)) -
F(x,zeta^2(t-h,x))}), we obtain that, possibly on a different set
$\wt\O$, with $P(\wt\O)=1$,
\begin{eqnarray*}
&   & G(x,\zeta_2(t-h,x))(B_t - B_{t-h}) \\
& = & G(x,\zeta_2(t,x))(B_t - B_{t-h})
      - (D_zG)(x,\zeta_2(t - h,x))G_2(x,\zeta_3(t - h,x))(B_t - B_{t-h})^2 \\
&   & + h^{3\alpha}R_\alpha(t-h,t,x) \\
& = & G(x,\zeta_2(t,x))(B_t - B_{t-h})
      - (D_zG)(x,\zeta_2(t,x))G_2(x,\zeta_3(t,x))(B_t - B_{t-h})^2 \\
&   & + h^{3\alpha}R_\alpha(t-h,t,x).
\end{eqnarray*}
Moreover, rewriting (\ref{F(x,zeta^2(t,x)) - F(x,zeta^2(t-h,x))}) as
$ F(x,\zeta_2(t-h,x))h   = F(x,\zeta_2(t,x))h +
h^{1+\alpha}R_\alpha(t-h,t,x) $, and noting from
(\ref{F(x,zeta^2(t,x)) - F(x,zeta^2(t-h,x))}) (recall definition
(\ref{Hz23})) that
\begin{eqnarray*}
&   & (D_zG)(x,\zeta_2(t-h,x))G_2(x,\zeta_3(t-h,x))(B_t - B_{t-h})^2 \\
& = & (D_zG)(x,\zeta_2(t,x))G_2(x,\zeta_3(t,x))(B_t - B_{t-h})^2
      + h^{3\alpha}R_\alpha(t-h,t,x),
\end{eqnarray*}
and that
\begin{eqnarray*}
&&(D_zG)(x,\zeta_2(t-h,x))G_2(x,\zeta_3(t-h,x))h \\
& = &(D_zG)(x,\zeta_2(t,x))G_2(x,\zeta_3(t,x))h
      + h^{1+\alpha}R_\alpha(t-h,t,x),
\end{eqnarray*}
we obtain from
(\ref{zeta(t,x) - zeta(t-h,x)}) that
\begin{eqnarray*}
&   & \zeta(t,x) - \zeta(t-h,x) \\
& = & \Big\{ F(x,\zeta_2(t,x)) - \frac{1}{2} (D_z
G)(x,\zeta_2(t,x))G_2(x,\zeta_3(t,x)) \Big\}h
      + G(x,\zeta_2(t,x))(B_t - B_{t-h}) \\
&   & - \frac{1}{2} (D_z G)(x,\zeta_2(t,x))G_2(x,\zeta_3(t,x))(B_t -
B_{t-h})^2
      + h^{3\alpha}R_\alpha(t-h,t,x).
\end{eqnarray*}
Finally, since $\alpha \in \big(\frac{1}{3},\frac{1}{2}\big)$ is
arbitrary, the proposition follows.
\qed

\subsection{Backward temporal-spatial expansion}

We now give the complete statement of the backward
temporal-spatial expansion.

\begin{prop}
\label{zeta(t-h,x+k) - zeta(t,x)}  Assume that $\zeta\in
C^{0,(3)}(\bF^B,[0,T],\hR^d)$ satisfying (T1)-(T3). Then, for all
$m\in\hN$ and $\alpha\in\left(\frac{1}{3},\frac{1}{2}\right)$ there
exists some subset $\widetilde\Omega \subset \Omega$ such that
$P(\widetilde\Omega)=1$, and that on $\widetilde\Omega$, for all $0
\leq t-h < t \leq T$, $x,k\in\hR^d$,
 \bea
 \label{dhzetats-}
 \zeta(t-h,x+k) - \zeta(t,x)& =& ah+b(B_{t-h} - B_t)+ \frac{c}{2}(B_{t-h} - B_t)^2+\lan p,k\ran
 +\frac{1}{2}\lan Xk,k\ran \nonumber\\
 &&+\lan q, k\ran (B_{t-h} - B_t)+ (h +
 |k|^2)^{3\alpha}R_\alpha(t-h,t,x,x+k),
 \eea
where
 \begin{eqnarray}
 a & = & -\left\{ F(x,\zeta_2(t,x)) - \frac{1}{2} (D_z
G)(x,\zeta_2(t,x))G_2(x,\zeta_3(t,x)) \right\};\nonumber\\
 b &=& G(x,\zeta_2(t,x)); \nonumber \\
 c & =  &  (D_z G)(x,\zeta_2(t,x))G_2(x,\zeta_3(t,x));\\
 p &=& (D_x\zeta)(t,x); \nonumber\\
 X &= &D_x^2\zeta(t,x);\nonumber \\
 q & =& (D_x G)(x,\zeta_2(t,x)) + (D_z G)(x,\zeta_2(t,x))D_x\zeta_2(t,x).\nonumber
 \end{eqnarray}
\end{prop}

{\it Proof.} As in the forward expansion case, we need only show
that
for all $m\in\hN$ and
$\alpha\in\left(\frac{1}{3},\frac{1}{2}\right)$ there exists some
subset $\widetilde\Omega \subset \Omega$ of full probability such
that on $\widetilde\Omega$, for all $0 \leq t-h < t \leq T$,
$x,k\in\hR^d$,
\begin{eqnarray}
\label{bckspexp}
 &&\zeta(t-h,x+k) - \zeta(t-h,x) \nonumber\\
& = & D_x\zeta(t,x)k + \frac{1}{2}\langle (D_x^2\zeta(t,x)k,k \rangle \\
&   & + \{(D_xG)(x,\zeta_2(t,x)) + (D_zG)(x,\zeta_2(t,x))D_x\zeta_2(t,x)\}k(B_{t-h} - B_t) \nonumber\\
&   & + (h + |k|^2)^{3\alpha} R_\alpha(t-h,t,x,x+k).\nonumber
\end{eqnarray}

From the usual Taylor expansion with the remainder in the Lagrange
form, we have
\begin{eqnarray*}
&   & \zeta(t-h,x+k) - \zeta(t-h,x) \\
& = & D_x\zeta(t-h,x)k + \frac{1}{2}\langle (D_x^2\zeta(t-h,x)k,k
\rangle + |k|^3R_\alpha(t-h,t,x,x+k),
\end{eqnarray*}
for $0 \leq t-h < t \leq T$, $x,k\in\hR^d$, $P$-a.s. Moreover,
from Proposition \ref{zeta(t-h,x) - zeta(t,x)} it follows that
\begin{eqnarray*}
&   & D_x\zeta(t,x) - D_x\zeta(t-h,x) \\
& = & D_x[G(x,\zeta_2(t,x))](B_t - B_{t-h}) + h^{2\alpha}R_\alpha(t-h,t,x) \\
& = & \{(D_xG)(x,\zeta_2(t,x)) +
(D_zG)(x,\zeta_2(t,x))D_x\zeta_2(t,x)\}(B_t - B_{t-h})
      + h^{2\alpha}R_\alpha(t-h,t,x)
\end{eqnarray*}
and
\[
  D_x^2\zeta(t,x) - D_x^2\zeta(t-h,x) = h^\alpha R_\alpha(t-h,t,x).
\]
Consequently,
\begin{eqnarray*}
&   & \zeta(t-h,x+k) - \zeta(t-h,x) \\
& = & D_x\zeta(t,x)k + \frac{1}{2}\langle (D_x^2\zeta(t,x)k,k \rangle \\
&   & + \{(D_xG)(x,\zeta_2(t,x)) + (D_zG)(x,\zeta_2(t,x))D_x\zeta_2(t,x)\}k(B_{t-h} - B_t) \\
&   & + |k|^{6\alpha}R_\alpha(t-h,t-h,x,k) + h^{2\alpha}k
R_\alpha(t-h,t,x)
      + h^\alpha\langle R_\alpha(t-h,t,x)k,k \rangle.
\end{eqnarray*}
Thus, by virtue of
\[
    |k|^{6\alpha} + h^{2\alpha}|k| + h^\alpha|k|^2 \leq C_m(h + |k|^2)^{3\alpha},
\]
for $h\in[0,T]$ and $|k| \leq m$, we derive (\ref{bckspexp}). This,
combined with (\ref{zeta(t,x) - zeta(t-h,x)}), leads us to the
backward temporal-spatial expansion (\ref{dhzetats-}).
\qed



Finally, by combining Proposition \ref{zeta(t+h,x+k) - zeta(t,x)}
and Proposition \ref{zeta(t-h,x+k) - zeta(t,x)}, we have completed
the proof of Theorem \ref{StoTaylor}.

\section{Application to Stochastic PDEs}
\setcounter{equation}{0}

\label{Notion of stochastic viscosity solution}

Having tried so hard to develop the various forms of stochastic
Taylor expansion, as an application in this section we shall try to
use it to study the {\it stochastic viscosity solutions} for fully
nonlinear SPDE, following the idea that we developed in our earlier
work \cite{Buckdahn_Ma_02}. In order not to over-complicate the computation
we shall consider the following simpler version of the fully nonlinear
SPDE (\ref{exspde}):
 \bea
 \label{SPDE}
  u(t,x) = u_0(x) + \int_0^t f(x,u,Du,D^2u) ds + \int_0^t g(x,Du)\circ dB_s,
  \q (t,x)\in[0,\infty)\times\hR^d.
 \eea
Compared to SPDE (\ref{exspde}), as well as Example \ref{egspde}, we see that
here the diffusion coefficient $g$ in (\ref{SPDE})
is independent of $u$. The general case could be treated in a similar way but with
more complicated expressions. Since our purpose here is to outline our idea of a
new definition of stochastic viscosity solution, without adding too much technical
complexity into this already lengthy paper, we shall leave the
study of the general case to a forthcoming paper.  We should also 
note that in
(\ref{SPDE}) we are using the Stratonovich integral instead of the It\^o integral
for the simplicity of the presentation, the following relation is worth noting:
 \bea
 \label{relat}
 g(x,Du)\circ dB_t=g(x,Du)dB_t+\frac12 D_zg(x,Du)D_x[g(s,Du)]dt.
 \eea
It is worth pointing out that even in this simplified form, the nonlinearity of the 
function $g$ on $Du$ already makes it difficult to apply the rough path approach of 
\cite{CFO} directly here.

To explain our idea of the definition of stochastic viscosity
solution, let us first apply Theorem \ref{StoTaylor}
to the regular solution $u$. Bearing the relation (\ref{relat}) in
mind we have:
\begin{eqnarray}
\label{u(t+h,x+k) - u(t,x)}
 u(t+h,x+k) - u(t,x) & = & a h +
b(B_{t+h}-B_t) + \frac{c}{2} (B_{t+h}-B_t)^2 + \langle p,k\rangle+
\frac{1}{2}\langle X k,k\rangle\nonumber
\\
&&+ \langle q,k\rangle (B_{t+h}-B_t)
+\left(|h|+|k|^2\right)^{3\alpha}R_{\alpha,m}(t,t+h,x,k),
\end{eqnarray}
where
 \bea
 \label{TaylorSPDE}
 a&=&  f(x,(u,Du,Du^2)(t,x)),\nonumber\\
 b&=&   g(x,Du(t,x))\\
 c&=&
 \lan D_zg(x,Du(t,x)),D_x[g(x,Du(t,x))]\ran, \nonumber\\
 p&=&  Du(t,x), \qq  X= D^2u(t,x),   \nonumber\\
 q&=& D_x[g((x,Du)(t,x)]=D_xg(x, Du)+
      D^2u(t,x) D_zg(x,Du(t,x)). \nonumber
 \eea
It is interesting to note that in this simple case  the terms
involving $D^2u$ can be written collectively as (suppressing
variables):
 \begin{eqnarray*}
&   & f(\cds, D^2u)+
\langle [D^2u] D_zg(x,Du)(B_{t+h}-B_t),k\rangle\\
&&+\frac{1}{2}\langle [D^2u] D_zg(x,Du), D_zg(x,Du)\rangle)
      (B_{t+h}-B_t)^2+ \frac{1}{2}\langle [D^2u] k,k\rangle\\
&=& f(\cds,D^2u)+\frac12\lan [D^2u] (D_zg(x,Du)(B_{t+h}-B_t)+k),
D_zg(x,Du)(B_{t+h}-B_t)+k\ran.
\end{eqnarray*}

Compared to the classical deterministic Taylor expansion, and the
SPDEs studied in \cite{Buckdahn_Ma_02} in which $g$ is independent
of $Du$, we can see that in a general case the terms involving
$D^2u$ becomes much more complicated, and
our previous method (via Doss-Sussmann) will face a fundamental
challenge, especially in the uniqueness proof. We therefore
will try to find a different approach to define the stochastic
viscosity solutions, using the stochastic characteristics introduced
by Kunita \cite{Kunita_90}, combined with our results on stochastic
Taylor expansions. This new method also reflects the basic ideas of
the works of Lions and Souganidis \cite{Lions_Souganidis_98_1},
\cite{Lions_Souganidis_98_2}, \cite{Lions_Souganidis_00}, and
\cite{Lions_Souganidis_02}, and in a sense includes our previous
work \cite{Buckdahn_Ma_02} as special case. For simplicity we shall
now assume all processes involved are real valued.


To motivate our definition of the stochastic viscosity solution let
us suppose in a first step that the coefficients of SPDE (\ref{Def.
u}) are sufficiently smooth and that this equation admits a regular
solution $u\in C^{0,\infty}(\bF^B,[0,T]\times \hR)$. Under this
assumption, we compare the solution $u$ with a smooth test field
$\varphi\in C^{0,\infty}(\bF^B,[0,T]\times \hR)$ defined as the
unique solution of the equation
\begin{eqnarray}
\label{Def. phi}
      d\varphi(t,x)
& = &  \theta(t,x)dt
      + g(x, D\varphi(t,x))\circ dB_t, \\
\nonumber
      \varphi(0,x)
& = & \varphi_0(x),
\end{eqnarray}
for $t\in[0,T]$, $x\in \hR^d$, where $\theta\in
C_{\ell,b}^{\infty}([0,T]\times \hR)$, $g\in C_{\ell,b}^{\infty}(\hR
\times \hR; \hR)$, and $\varphi_0\in C_p^{\infty} (\hR)$.

Let us now fix an $\cF^B$-measurable $[0,T]$-valued random variable
$\tau$ and an $\cF^B$-measurable $\hR^d$-valued random variable
$\xi$. We say that $u-\varphi$ achieves a {\it local left-maximum}
in $(\tau,\xi)$ if for almost all $\omega\in\{\tau<T\}$ there is
some $\rho>0$ (which may depend on $\omega$) such that
\[
  (u-\varphi)(\omega,t,x)\leq (u-\varphi)(\omega,\tau(\omega),\xi(\omega)),
\]
for all $(t,x)\in [0,T]\times \hR^d$ with $t\in((\t(\omega)-\rho)^+,
\t(\o)]$ and $|x-\xi(\omega)|\leq \rho$.

Following the approach by Kunita (see Theorem 6.1.2 in
\cite{Kunita_90}), we formally introduce the following stochastic
characteristics. For the moment let us assume that the random time
$\t$ is actually deterministic to avoid further complications (keep
in mind, however, that the Taylor expansion will hold even for the
arbitrary random time $\t$!).
 \bea
  \label{stochar}
 \f_t(x,z) & =& \dis x - \int_\tau^t D_zg((\f_s,\chi_s)(x,z))\circ dB_s,
 \nonumber\\
 \eta_t(x,y,z) &= &\dis y +\int_\tau^t \{g((\f_s,\chi_s)(x,z))
 - \chi_s(x,z)D_z g((\f_s,\chi_s)(x,z))\}\circ dB_s, \\
 \chi_t(x,z) & = &\dis z + \int_\tau^t D_x g((\f_s,\chi_s)(x,z))
 \circ dB_s, \nonumber
 \eea
where $t\in[0,\tau]$, $(x,y,z)\in\hR^3$, and we hope to be able to
define a transformation $\psi(t,x)$ by
\begin{equation}
  \label{Formal def. psi}
  \varphi(t,\f_t(x,D\psi(t,x))) =
  \eta_t(x,\psi(t,x),D\psi(t,x)).
\end{equation}
For notational simplicity let us now denote:
 \bea
 \label{fhk}
 l(x,z)=-D_zg(x,z);\q
 h(x,z)= g(x,z)  - zD_z g(x,z);\q
 k(x,z)=D_x g(x,z).
 \eea
\begin{rem} {\rm We note that if the function $g$ is linear in $Du$, then the situation will become
much simpler, and the arguments below would become straightforward. To be more precise, let us
consider the following two cases:

 (i)  $g(x,Du)=\lan H, Du\ran$, where $H$ is a constant vector. Since $D_xg=D_xV(x)$ and $D_zg=H$, 
the Taylor expansion (\ref{TaylorSPDE}) is drastically simplified. Also, the characteristics
(\ref{stochar}) becomes almost triviel: $\chi\equiv z$, $\eta\equiv y$, and $\f\equiv x-H(t-\t)$.

\ms

(ii) $g(x,Du)=V(x)Du$ (see \cite{CFO}). In this case one has $D_xg(x,z)=D_xV(x)z$, $D_zg(x,z)=V(x)$, 
thus the Taylor expansion (\ref{TaylorSPDE}) will also become much simpler. Furthermore, the characteristics (\ref{stochar}) now become ``disentangled" SDEs:
 \beaa
 \phi_t(x)&=&x-\int_\tau^t V(\phi_s(x))\circ dB_s,\\
 \eta_t(x,y,z)&\equiv&y,\\
 \chi_t(x,z)&=&z+\int_\t^t D_xV(\phi_s(x))\chi_s(x,z)\circ dB_s,
 \eeaa
and (\ref{Formal def. psi}) takes the special form: $\psi(t,x)=\varphi(t,\phi_t(x))$. In other words, the transformation $\psi$ is globally well-defined  by an easy and explicit expression. \qed
}
\end{rem}
In what follows we shall denote
$\Th\dfnn (l,h,k)$. Also, for any function $\g=\g(x,y,z)$ we denote
$\nabla \g= (D_x\g,D_y\g,D_z\g)^T$. Then we can write, for example,
the first equation in (\ref{stochar}) in the It\^o integral from:
 \bea
 \label{stochar2}
 \f_t(x,z)&=&x+ \int_\t^tl((\f_s,\chi_s)(x,z))\circ dB_s\\
  &=&x+\int_\t^t\frac12 \lan\nabla l,\Th\ran((\f_s,\chi_s)(x,z))ds
  +\int_\t^tl((\f_s,\chi_s)(x,z))  dB_s. \nonumber
 \eea

Now, treating the random fields $\eta$ and $\chi$ in (\ref{stochar})
the same way, and applying the stochastic (backward temporal) Taylor
expansion to $(\f,\eta,\chi)$
we have, for all $\a\in(\frac13,\frac12)$,
 \bea
 \label{stocharexpphi}
 \left\{\begin{array}{lll}
 \dis \f_t(x,z)=x+l(x,z)(B_t-B_\t)+\frac12\lan \nabla l, \Th\ran(x,z)(B_t-B_\t)^2\
 +|t-\t|^{3\a} R_{\a,m},\\
 \dis \eta_t(x,y,z) = y + h(x,z)(B_t-B_\t)+\frac12\lan\nabla
 h,\Th\ran(x,z)(B_t-B_\t)^2+|t-\t|^{3\a} R_{\a,m}, \\
 \dis \chi_t(x,z)  =  z + k(x,z)(B_t-B_\t)+\frac12\lan\nabla
 k,\Th\ran(x,z)(B_t-B_\t)^2+|t-\t|^{3\a} R_{\a,m}.
 \ea\right.
 \eea

On the other hand, writing (\ref{Def. phi}) in an It\^o integral
form we have
 \bea
 \label{varphieq}
 \varphi(t,x)&=&\varphi_0(x)+\int_\t^t \{\th(s,x)+\frac12
 D_zg(x,D\varphi(s,x))D_x[g(x,D\varphi(s,x))]\}ds\nonumber\\
 &&+\int_\t^t g(x,D\varphi(s,x))dB_s  \\
 &\dfnn& \varphi_0(x)+\int_\t^t F(s,x)ds+\int_\t^tG(x,
 \zeta_2(s,x))dB_s.\nonumber
 \eea
where $F$ and $G$ are defined in an obvious way, and with
 $$ \zeta_1(t,x)\dfnn\varphi(t,x), \q \zeta_2(t,x)\dfnn D\varphi(t,x),\q
 \zeta_3(t,x)=(D\varphi, D^2\varphi)(t,x).
 $$
Note that differentiating the both sides of (\ref{varphieq}) we have
 \bea
  dD\varphi(t,x)=
  D_xF(t,x)dt+\{D_xg(x,D\varphi(t,x))+D_zg(x,D\varphi(t,x))D^2_{x}
  \varphi(t,x))\}dB_t.
 \eea
 Applying the backward temporal Taylor expansion
(Theorem 2.3) again to the random field $\varphi$ around any point
$(\t,\vsi)$, with
 $$
 F_2(t,x)\dfnn D_xF(t,x);\qq G_2(x,z_1,z_2)\dfnn
 D_xg(x,z_1)+D_zg(x,z_1)z_2,
 $$
we obtain, after some simple cancelations, that
 \bea
 \label{varphiexp}
 \varphi(t,x)&=&\varphi(\t,\vsi)+
 \theta(\t,\vsi)(t-\t)+g(\vsi,D\varphi(\t,\vsi))(B_t-B_\t)\nonumber\\
 &&+\frac12\Big\{[D_z g
 D_xg](\vsi,D\varphi(\t,\vsi))+(D_zg(\vsi, D\varphi(\t,\vsi)))^2D^2_{x}\varphi(\t,\vsi)
 \Big\}(B_t-B_\t)^2\nonumber\\
 &&+D\varphi(\t,\vsi)(x-\vsi)+
 D_x[g(\cd,D\varphi(\cd,\cd))](\t,\vsi)(x-\vsi)(B_t-B_\t)\\
 && +\frac12 D^2_x\varphi(\t,\vsi)(x-\vsi)^2+(|t-\t|+|x-\vsi|^2)^{3\a}R_{\a,m}.
\nonumber
 \eea
We note that the above holds for all $(\t,\vsi)$ and all
$\o\in\widetilde \O_{\a,m}$.

Let us now define the desired random field $\psi$. We should note
that our main purpose here is to find such a transformation so as to
eliminate the stochastic integral. In other words, we shall look for
such $\psi$ that has the following first order Taylor expansion:
 \beq
 \label{psiform}
 \left\{\begin{array}{lll}
 \dis \psi(t,x)=\varphi(\t,x)+ \pa^-_t\psi(\t,x)(t-\t)
 +|t-\t|^{3\a}R_{\a,m};\\
\dis D\psi(t,x)=D\varphi(\t,x)+\pa^-_tD\psi(\t,x)(t-\t)
 +|t-\t|^{3\a}R_{\a,m}.
 \end{array}
 \right.
 \eeq
Here $\pa^-_t$ denotes the left partial derivative with respect to
$t$. To this end, we note that (\ref{stochar}) and (\ref{Formal def.
psi}) imply that for $t=\t$, one has
 \bea
 \label{phieqpsi}
 \varphi(\t,x)=\varphi(\t,\f_\t(x,D\psi(\t,x)))=\eta_\t
(x,\psi(\t,x),D\psi(\t,x))=\psi(\t,x), ~~x\in\hR,
 \eea
and hence $D\varphi(\t,x)=D\psi(\t,x)$ holds for all $x\in\hR$ as
well. Next, we look at the Taylor expansion for both $\f$ and
$\eta$. Recalling (\ref{fhk}) and the Taylor expansions
(\ref{stocharexpphi}). We shall first take $\bx= (x,\psi(t,x),
D\psi(t,x)))$, and then  replacing $x$ by $\f_t(x,D\psi(t,x))$ and
$\vsi$ by $x$ in (\ref{varphiexp}). It should be noted that after
these substitutions the remainder will  look like
$R_{\alpha,m}(t,\tau,x,\phi_t(x,D\psi(t,x)))$, a slightly more
complicated form than the original ones. To make sure the accuracy
of the expansion in what follows we shall denote
$$\overline{R}_{\alpha,m,loc}=\sup_{t,s\in[0,T];~
x\in\overline{B}_m(0)}|R_{\alpha,m}(t,s,x,\phi_t(x,D\psi(t,x)))|.$$
Then it is not hard to show that there exists an increasing sequence
$\{\Omega_\ell,\ell\ge 1\}\subset {\cal F}^B$ with
$\lim_{\ell\rightarrow +\infty}P(\Omega_\ell)= 1$,   such that,
$$\overline{R}_{\alpha,m,loc}I_{\Omega_\ell}\in L^{\infty}_-(\O, \cF,
P),\quad\mbox{for all~~}\, \ell\ge 1.$$

Keeping such a modification in mind, we can now proceed to write
down the Taylor expansion:
 \bea
 \label{varphiexp2}
&&
\varphi(t,\f_t(x,D\psi(t,x)))\nonumber\\
&=&\varphi(\t,x)+
 \theta(\t,x)(t-\t)+g(x,D\varphi(\t,x))(B_t-B_\t)\nonumber\\
 &&+\frac12\Big\{[D_z g D_xg](x,D\varphi(\t,x))+(D_zg(x, D\varphi(\t,x))^2D^2_x
 \varphi(\t,x)\Big\}(B_t-B_\t)^2\nonumber\\
 &&-D\varphi(\t,x)D_zg(x,D\varphi(\t,x))(B_t-B_\t)\nonumber\\
 &&-D\varphi(\t,x)D^2_zg(x,D\varphi(\t,x))\pa_t^-D\psi(\t,x)(t-\t)(B_t-B_\t)\nonumber\\
 &&+\frac12 D\varphi(\t,x)\lan \nabla l,\Th\ran(x,
 D\varphi(\t,x))(B_t-B_\t)^2\\
 &&-\frac12 D^2_x\varphi(\t,x) (D_zg(x, D\varphi(\t,x)))^2(B_t-B_\t)^2 \nonumber\\
 &&-D_xg(x,D\varphi(\t,x))D_zg(x,D\varphi(\t,x))(B_t-B_\t)^2
 +|t-\t|^{3\a}\overline{R}_{\a,m,loc}\nonumber\\
 &=&\varphi(\t,x)+ \theta(\t,x)(t-\t)+g(x,D\varphi(\t,x))(B_t-B_\t)\nonumber\\
 &&-\frac12[D_z g D_xg](x,D\varphi(\t,x))(B_t-B_\t)^2-D\varphi(\t,x)D_zg(x,
 D\varphi(\t,x))(B_t-B_\t)\nonumber\\
 &&+\frac12 D\varphi(\t,x)\lan \nabla l,\Th\ran(x,
 D\varphi(\t,x))(B_t-B_\t)^2 \nonumber\\
 &&-D\varphi(\t,x)D^2_zg(x,D\varphi(\t,x))\pa_t^-D\psi(\t,x)(t-\t)(B_t-B_\t)+|t-\t|^{3\a}
 \overline{R}_{\a,m,loc}.
 \nonumber
 \eea
On the other hand, from
(\ref{stochar}) we see that
 \beaa
 \psi(t,x)&=&\eta_t(x,\psi(t,x),D\psi(t,x))-h(x,
 D\psi(t,x))(B_t-B_\t)\\
 &&-\frac12 \lan \nabla h,\Th\ran(x,D\psi(t,x))
 (B_t-B_\t)^2+|t-\t|^{3\a}\overline{R}_{\a,m,loc}.
 \eeaa
Using the form (\ref{psiform}) and the smoothness assumptions on all
the coefficients, and then noting (\ref{phieqpsi}) one can easily
rewrite above as
 \beaa
 \psi(t,x)&=&\eta_t(x,\psi(t,x),D\psi(t,x))-h(x,
 D\varphi(\t,x))(B_t-B_\t)\\
 &&-D_zh(x,D\varphi(\t,x))\pa^-_tD\psi(\t,x)(t-\t)(B_t-B_\t)\\
 &&-\frac12 \lan \nabla h,\Th\ran(x,D\varphi(\t,x))
 (B_t-B_\t)^2+|t-\t|^{3\a}\overline{R}_{\a,m,loc}.
 \eeaa
Now recalling the relation (\ref{Formal def. psi}), the fact
(\ref{varphiexp2}), and the definition of $h$ (\ref{fhk}), we obtain
that
 \bea
 \label{psidef}
 \psi(t,x)&=&\eta_t(x,\psi(t,x),D\psi(t,x))-h(x,
 D\varphi(\t,x))(B_t-B_\t)\nonumber\\
 &&-\frac12 \lan \nabla h,\Th\ran(x,D\varphi(\t,x))
 (B_t-B_\t)^2+|t-\t|^{3\a}R_{\a,m}\nonumber\\
 &=&\varphi(\t,x)+ \theta(\t,x)(t-\t)+g(x,D\varphi(\t,x))(B_t-B_\t)\nonumber\\
 &&-\frac12[D_z g D_xg](x,D\varphi(\t,x))(B_t-B_\t)^2-D\varphi(\t,x)D_zg(x,
 D\varphi(\t,x))(B_t-B_\t)\nonumber\\
 &&+\frac12 D\varphi(\t,x)\lan \nabla l,\Th\ran(x,
 D\varphi(\t,x))(B_t-B_\t)^2\nonumber\\
 &&-D\varphi(\t,x)D^2_xg(x,D\varphi(\t,x))\pa_t^-D\psi(\t,x)(t-\t)(B_t-B_\t)\\
 &&-[g(x, D\varphi(\t,x))-D\varphi(\t,x)D_zg(x,D\varphi(\t,x))](B_t-B_\t)\nonumber\\
 &&-D_zh(x,D\varphi(\t,x))\pa^-_tD\psi(\t,x)(t-\t)(B_t-B_\t)\nonumber\\
 &&-\frac12 \lan \nabla h,\Th\ran(x,D\varphi(\t,x))
 (B_t-B_\t)^2+|t-\t|^{3\a}\overline{R}_{\a,m,loc}\nonumber\\
 &=&\varphi(\t,x)+ \theta(\t,x)(t-\t)-\frac12[D_z g D_xg](x,D\varphi(\t,x))
 (B_t-B_\t)^2\nonumber\\
 &&+\frac12 D\varphi(\t,x)\lan \nabla l,\Th\ran(x,
 D\varphi(\t,x))(B_t-B_\t)^2\nonumber\\
 &&-D\varphi(\t,x)D^2_xg(x,D\varphi(\t,x))\pa_t^-D\psi(\t,x)(t-\t)(B_t-B_\t)\nonumber\\
 &&-D_zh(x,D\varphi(\t,x))\pa^-_tD\psi(\t,x)(t-\t)(B_t-B_\t)\nonumber\\
  &&-\frac12 \lan \nabla h,\Th\ran(x,D\varphi(\t,x))
 (B_t-B_\t)^2+|t-\t|^{3\a}\overline{R}_{\a,m,loc}.\nonumber
 \eea
Since $\nabla l=(-D_{xz}g,0,- D_z^2g)$ and $\nabla
h=(D_xg-zD_{xz}g,0,
 -zD_z^2g)$, one can check that
 \beaa
 -\frac12[D_z g D_xg](x,D\varphi(\t,x))
 +\frac12 D\varphi(\t,x)\lan \nabla l,\Th\ran(x,
 D\varphi(\t,x))-\frac12 \lan \nabla h,\Th\ran(x,D\varphi(\t,x))
 =0,
 \eeaa
and
 \beaa
 &&-D\varphi(\tau,x)D^2_zg(x,D\varphi(\tau,x))\partial^-_tD\psi
 (\tau,x))(t-\tau)(B_t-B_\tau)\\
 &&\qq\qq -D_zh(x,D\varphi(\tau,x))\partial_t^-D\psi(\tau,x)(t-\tau)
 (B_t-B_\tau)=0.
 \eeaa
In other words, (\ref{psidef}) leads to that
 \bea
 \label{defpsi}
 \psi(t,x)=\varphi(\t,x)+\th(\t,x)(t-\t)+|t-\t|^{3\a}\overline{R}_{\a,m,loc}.
 \eea
Clearly, (\ref{defpsi}) indicates that
 $$ \pa^-_t\psi(\t,x)=\th(\t,x),\qq \pa^-_tD\psi(\t,x)=D\th(\t,x).
 $$


We can now define the notion of stochastic viscosity solution for
(\ref{Def. u}).

\begin{defn}
\label{DefSVS} A random field $u\in C([0,T]\times\hR^d)$ is called a
stochastic viscosity subsolution (resp.~supersolution) of (\ref{Def.
u}), if $u(0,x) \leq u_0(x)$ ($u(0,x) \geq u_0(x)$) for all
$x\in\hR^d$ and if for any $\tau\in L^0(\cF^B,[0,T])$, $\xi\in
L^0(\cF^B,\hR^d)$, and any (not necessarily adapted) random field
$\varphi\in C^{0,2}([0,T]\times\hR^d)$ having the expansion
(\ref{varphiexp}),  it holds that
\begin{eqnarray}
\nonumber
    \th(\tau,\xi)
& \leq (\geq) &
f(\xi,u(\tau,\xi),D\varphi(\tau,\xi),D^2_x\varphi(\t,\xi)).
\end{eqnarray}
 $P$-a.s. on the subset of $\O$ on which $u-\varphi \le (\ge) (u-\varphi)(\t,\xi)=0$
 at a left neighborhood of $(\t,\xi)$. A random field $u\in
C([0,T]\times\hR^d)$ is called a stochastic viscosity solution of
(\ref{Def. u}), if it is both a stochastic viscosity subsolution and
a supersolution.
\end{defn}

\begin{rem}{\rm We would like to note that in Definition \ref{DefSVS} the
(possibly anticipating) test functions $\varphi$ can be directly
defined by the expansion (\ref{varphiexp}) with $\varphi(\tau,\xi)$
being replaced by $u(\tau,\xi)$ and
$(\theta(\tau,\xi),D_x\varphi(\tau,\xi),D^2_x\varphi(\tau, \xi))$
being replaced by a triplet of random variables $(\beta, p, A)$. The
main advantage here is that the random field $\varphi$ is now
defined {\it globally}, overcoming the essential difficulties in the
theory of stochastic viscosity solution thus far, and will
significantly facilitate the uniqueness proof. Moreover, on the
subset of $\Omega$ where $u-\varphi$ achieves the local left-maximum
(resp., left-minimum) at the point $(\tau,\xi)$, the triplet
$(\beta, p, A)$ can be considered as a stochastic sub-
(resp.super-)jet, as it was traditionally done. These issues will be
further explored in our forthcoming publications. } \qed
\end{rem}

In the rest of this section we shall verify
that a regular solution must be a stochastic viscosity solution in
the sense of Definition \ref{DefSVS}, which will
provide a justification for our new definition.
To this end, let us assume that the coefficient $f$ is {\it proper}
in the following sense

\ms

 {\bf (H1)} {\it The function $F(t,u,p,X)\dfnn -f(x,u,p,X)$ is ``{\it
 degenerate elliptic}". That is, $f$ is continuous in all
variables, and is {\it non-decreasing} in the variable $X$. }

\ms

Assume that $u$ is a regular solution to (\ref{Def. u}), then we
have
 \bea
 \label{e'}
 u(t,x) = u_0(x)+\int_0^t F(s,x)ds+ \int_0^tg(x,Du(s,x)) dB_s,
 \qq t\ge 0,
 \eea
where
 \beaa
 F(t,x)=f(x,(u,Du,D^2u)(t,x))+\frac{1}{2}D_zg(x,Du(t,x))D_x\left[g(\cd,Du(t,\cd))
 \right](t,x).
 \eeaa
For any given pair of random variables $(\tau,\xi)$ and an arbitrary
test field $\varphi$ such that $u-\varphi$ attains a local
left-maximum at $(\t,\xi)$ on a subset of $\O$ with positive
probability, let $\psi$ be the process associated to $\varphi$ by
(\ref{defpsi}) (where the remainder can be chosen to be zero!), and
$\phi$ be defined by (\ref{stocharexpphi}).

We first apply the Taylor expansion on $u$ at point $(\t,x)$, and
evaluated at $(t,\f_t(x, D\psi(t,x)))$ and then use the expansion of
$\phi_t(t,D\psi(t,x))$ (recall (\ref{stocharexpphi})) to get
 \bea
 \label{clas}
 && u(t,\phi_t(x,D\psi(t,x))) \nonumber\\
 & =& u(\tau,x)+f(x,(u,Du,D^2u)(\tau,x))(t-\tau)+g(x,Du(\tau,x))(B_t-B_\tau)\nonumber\\
 &&+\frac{1}{2}\left\{
 D_zg(x,Du(\tau,x))D_xg(x,Du(\tau,x))
 +D_zg(x,Du(\tau,x))^2D_{xx}^2u(t,x)\right\}(B_t-B_\tau)^2\nonumber \\
 &&+D_xu(\tau,x)(\phi_t(x,D\psi(t,x))-x)+\frac{1}{2}D_{xx}^2u(\tau,x)(\phi_t(x,D\psi(t,x))
 -x)^2\nonumber\\
 && +D_x[g(x,Du(\tau,x))](B_t-B_\tau)(\phi_t(x,D\psi(t,x))-x)
 +|t-\tau|^{3\alpha}R_{\a, m}(\tau,t,x,\f_t(x,D\psi(t,x)))
 \nonumber\\
 &=& u(\tau,x)+f(x,(u,Du,D^2u)(\tau,x)(t-\tau)+
 g(x,Du(\tau,x))(B_t-B_\tau)\nonumber\\
 && +\frac12\Big\{[D_z g D_xg](x,Du(\tau,x))+(D_zg(x,
 Du(\tau,x))^2D^2_xu(\tau,x)\Big\}(B_t-B_\tau)^2\\
 &&-Du(\tau,x)D_zg(x,D\varphi(\tau,x))(B_t-B_\tau)\nonumber\\
 && -Du(\tau,x)D^2_zg(x,D\varphi(\tau,x))\partial^-_tD\psi
 (\tau,x))(t-\tau)(B_t-B_\tau)\nonumber\\
 &&+\frac12 Du(\tau,x)\langle \nabla f,\Theta\rangle(x,
 D\varphi(\tau,x))(B_t-B_\tau)^2\nonumber\\
 && +\frac12 D^2_xu(\tau,x) (D_zg(x,
 D\varphi(\tau,x)))^2(B_t-B_\tau)^2\nonumber\\
 &&-(D_xg)(x,Du(\tau,x))D_zg(x,D\varphi(\tau,x))(B_t-B_\tau)^2\nonumber\\
 &&-(D_zg)(x,Du(\tau,x))D_zg(x,D\varphi(\tau,x))D^2_xu(\tau,x)(B_t-B_\tau)^2
 +|t-\tau|^{3\alpha}\overline{R}_{\alpha,m,loc}.\nonumber
 \eea
We will now replace $x$ by the random point $\xi$. To specify the
``left-neighhood" required for the viscosity property, we also
define for any $\rho>0$ the following subset of $\O$:
$$\Gamma_{\tau,\xi}^{\varphi,\rho}\dfnn \{\omega: (u-\varphi)(\o,t,x)\le
(u-\varphi)(\omega,\tau(\omega),\xi(\omega)),\mbox{~~$(t-\t(\o))^-<\rho$
and $x\in B_\rho(\xi(\o))$}\},
$$
To wit, $\G^{\varphi,\rho}_{\t,\xi}$ is the subset of $\O$ on which
$u-\varphi$ attains a {\it left local maximum} at $(\t,\xi)$. Now
setting $x=\xi$ in (\ref{clas}), and noting that
$Du(\tau,\xi)=D\varphi(\tau,\xi)$ and
$u(\tau,\xi)=\varphi(\tau,\xi)$, we have
 \beaa
 &&u(t,\phi_t(\xi,D\psi(t,\xi)))\\
 &=&u(\tau,\xi)+ f(\xi,(u,Du,D^2u)(\tau,\xi))(t-\tau)+g(\xi,
 Du(\tau,\xi))(B_t-B_\tau)\\
 && -\frac12[D_z g D_xg](\xi,Du(\tau,\xi))(B_t-B_\tau)^2
 -Du(\tau,x)D_zg(\xi,Du(\tau,\xi))(B_t-B_\tau)\\
 &&+\frac12 Du(\tau,\xi)\langle \nabla f,\Theta\rangle(x,
 Du(\tau,\xi))(B_t-B_\tau)^2\\
 &&-Du(\tau,\xi)D^2_zg(\xi,Du(\tau,\xi))\partial^-_tD\psi
 ( \tau,\xi))(t-\tau)(B_t-B_\tau)+|t-\tau|^{3\alpha}R_{\alpha,m,loc},
 \eeaa
almost surely  on $\Gamma_{\tau,\xi}^{\varphi,\rho}$. Consequently,
from (\ref{varphiexp2}) we obtain that, $P$-a.s. on
$\Gamma_{\tau,\xi}^\varphi$,
 \beaa
 0&\ge&\big(u(t,\phi_t(\xi,D\psi (t,\xi)))
 -u(\tau,\xi)\big)-\big(\varphi(t,\phi_t(\xi,D\psi (t,\xi)))
 -\varphi(\tau,\xi)\big)\\
 &=&\{f(\xi,(u,Du,D^2u)(\tau,\xi))-\theta(\tau,\xi)\}(t-\tau)
 +|t-\tau|^{3\alpha}\overline{R}_{\alpha,m,loc}. \
 \eeaa
Since $t\le\tau$ on $\Gamma_{\tau,\xi}^\varphi$, we deduce that
 $$f(\xi,(u,Du,D^2u)(\tau,\xi))-\theta(\tau,\xi)\ge 0, \q \mbox{$P$-a.s. on } \Gamma_{\tau,
 \xi}^\varphi.$$
Finally, since the mapping $X\mapsto f(x,u,p,X)$ is non-decreasing,
thanks to (H1),
and since $D^2\varphi(\tau,\xi)\ge D^2u(\tau,\xi)$ on
$\Gamma_{\tau,\xi}^\varphi$ we have
 $$f(\xi,(\varphi,D\varphi,D^2\varphi)(\tau,\xi))\ge
 \theta(\tau,\xi)\q \mbox{$P$-a.s. on } \Gamma_{\tau,
 \xi}^\varphi.$$
This proves that the classical solution $u$ is a stochastic
viscosity subsolution. That $u$ is also a supersolution can be
proved using a similar argument. Therefore, $u$ is a stochastic
viscosity solution.



\end{document}